\documentclass[12pt,reqno]{amsart}

\usepackage[T1]{fontenc}
\usepackage{lmodern}
\usepackage{amsmath,amssymb,amsthm,mathtools}
\usepackage[margin=1in]{geometry}
\usepackage{microtype,booktabs,array,enumitem}
\usepackage{graphicx,xcolor}
\usepackage[unicode,colorlinks=true,linkcolor=blue!50!black,citecolor=blue!50!black,urlcolor=blue!50!black]{hyperref}
\hypersetup{pdftitle={Split Hamiltonian structure for Vlasov--Maxwell equations and cylindrical pinch equilibria},pdfauthor={Justin Holmer and Katherine Zhiyuan Zhang}}
\numberwithin{equation}{section}
\newtheorem{theorem}{Theorem}[section]
\newtheorem{proposition}[theorem]{Proposition}
\newtheorem{lemma}[theorem]{Lemma}
\newtheorem{corollary}[theorem]{Corollary}
\theoremstyle{definition}
\newtheorem{assumption}[theorem]{Assumption}

\theoremstyle{remark}
\newtheorem{remark}[theorem]{Remark}
\DeclareMathOperator{\curl}{curl}
\DeclareMathOperator{\diver}{div}
\DeclareMathOperator{\Ran}{Ran}
\DeclareMathOperator{\Ker}{ker}
\DeclareMathOperator{\spanop}{span}
\DeclareMathOperator{\esssup}{ess\,sup}
\DeclareMathOperator{\sech}{sech}
\newcommand{\R}{\mathbb R}

\newcommand{\T}{\mathbb T}
\newcommand{\cK}{\mathcal K}
\newcommand{\cI}{\mathcal I}
\newcommand{\cE}{\mathcal E}
\newcommand{\cV}{\mathcal V}
\newcommand{\cJ}{\mathcal J}
\newcommand{\cL}{\mathbf L}
\newcommand{\cA}{\mathcal A}
\newcommand{\cB}{\mathcal B}
\newcommand{\cX}{\mathcal X}
\newcommand{\cY}{\mathcal Y}
\newcommand{\cZ}{\mathcal Z}
\newcommand{\cF}{\mathcal F}
\newcommand{\cT}{\mathcal T}
\newcommand{\VV}{\parallel}
\newcommand{\HH}{\perp}
\newcommand{\ev}{\mathrm{ev}}
\newcommand{\od}{\mathrm{od}}
\newcommand{\pre}{\mathrm{pre}}
\newcommand{\rev}{\mathrm{rev}}
\newcommand{\acc}{\mathrm{acc}}
\newcommand{\kin}{\mathrm{kin}}
\newcommand{\ip}[2]{\langle #1,#2\rangle}
\newcommand{\norm}[1]{\lVert #1\rVert}
\newcommand{\nm}{n^-}
\newcommand{\cl}{\overline}
\newcommand{\dd}{\,d}
\newcommand{\eeta}{e_\eta}
\newcommand{\esig}{e_\sigma}
\newcommand{\eR}{e_r}
\newcommand{\eTh}{e_\theta}
\newcommand{\eZ}{e_z}
\allowdisplaybreaks[2]
\title[Split Hamiltonian structure and pinch equilibria]{Split Hamiltonian structure for Vlasov--Maxwell equations and cylindrical pinch equilibria}
\author{Justin Holmer}
\address{Department of Mathematics, Brown University}
\author{Katherine Zhiyuan Zhang}
\address{Department of Mathematics, Northeastern University}
\dedicatory{Dedicated to Walter A. Strauss}
\date{November 19, 2026}

\begin{document}

\begin{abstract}
We study the linearized Vlasov--Maxwell equations for cylindrical pinch equilibria, including the classical Bennett pinch, Gaussian equilibria from our earlier work, and a rigidly rotating $\theta$-pinch. Conservation of particle momentum along a symmetry brings the energy principle of Lin and Strauss into these geometries. We work on a graph space for the symmetry component of the vector potential. For the Bennett $z$-pinch the bare potential form has one negative direction, but the particle-number constraint removes it and the entire $z$-independent spectrum lies on the imaginary axis. For a $z$-pinch with variable axial drift, every $z$-independent growing mode is radial. We prove axisymmetric spectral stability for a Gaussian $\theta$-pinch in an explicit range of axis slopes and for a center-peaked, rigidly rotating $\theta$-pinch. In a screw-pinch the helical twist couples the magnetic flux to the particle and electric constraints. The radial equations still separate, and all sufficiently high helical harmonics have imaginary spectrum. For two Gaussian screw-pinch equilibria we verify spectral stability in the full helical symmetry class.
\end{abstract}
\maketitle
\tableofcontents

\section{Introduction}\label{sec:intro}

Can a small disturbance of a plasma column grow exponentially? In the Bennett $z$-pinch an axial current produces the magnetic field that confines the particles \cite{Bennett1934}. Before the particle constraints are imposed, the quadratic form in the magnetic potential has one negative direction. Its negative term comes from the equilibrium current response. A growing mode, however, preserves the particle number of each species to first order. We show that this constraint alone removes the negative direction. Accounting for the constraints that particle orbits impose on a field variation is central to the energy principle of Zhiwu Lin and Walter Strauss \cite{LinStrauss2007,LinStrauss2008}.

\begin{figure}[tb]
\centering
\includegraphics[width=.87\linewidth]{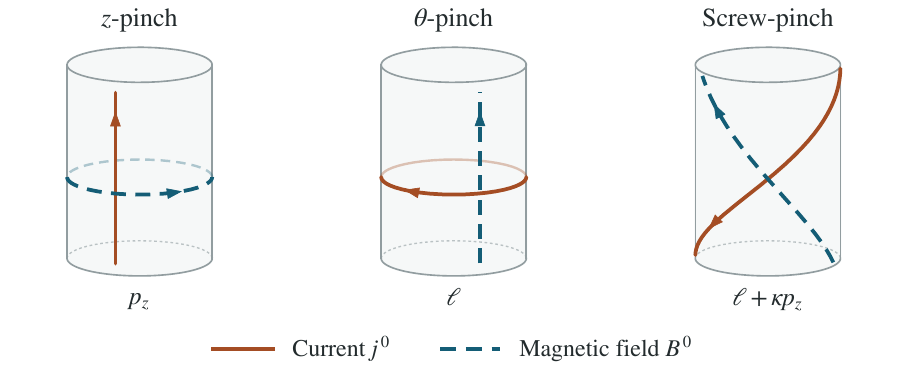}
\caption{Current and magnetic field on a cylinder of fixed radius. Below each column is the conserved momentum for the symmetry shown. The solid curve follows the current and the dashed curve follows the magnetic field. In the screw-pinch panel the current follows a helix of pitch parameter $\kappa$, and the common-pitch equilibrium field is perpendicular to it.}
\label{fig:pinch-geometry}
\end{figure}

An earlier theorem of Glassey and Strauss explains our choice of background. For a finite-energy plasma on $\R^3$ without external fields, the energy in each fixed ball tends to zero along a sequence of times \cite[Theorem~2]{GlasseyStrauss1984}. A stationary plasma with this property must be trivial. In earlier work we proved the corresponding stationary nonexistence result on $\R^2\times\T$ and studied pinch equilibria of infinite background energy \cite[Proposition~A.1 and Theorems~1.1--1.4]{HolmerZhang2026}. Here we study finite-energy disturbances with the symmetries in Figure~\ref{fig:pinch-geometry}. On the whole cylinder we take the symmetry component $a$ of the vector potential in $D(C_e)$, the graph space of its curl block, with norm $(\|a\|^2+\|C_ea\|^2)^{1/2}$. All spectral assertions refer to this realization. For every physical finite-energy normal mode with spectral parameter $\lambda\ne0$ in the kinetic and Maxwell domains specified in Section~\ref{sec:framework}, Faraday's law gives such a potential. By spectral stability we mean that the spectrum of the generator is contained in $i\R$.

For relativistic equilibria decreasing in particle energy, Lin and Strauss reduced the growing-mode problem in periodic and axisymmetric geometries to a self-adjoint operator on a single field potential. In their notation it is
\[
 \mathcal L^0=(\mathcal B^0)^*(\mathcal A_1^0)^{-1}\mathcal B^0+\mathcal A_2^0 .
\]
The inverse solves the electrostatic constraint after fixing the gauge. The operators retain the projection onto functions constant along equilibrium particle orbits. Nonnegativity gives spectral stability \cite{LinStrauss2007}. Their sharp criterion counts independent growing modes by the negative eigenvalues of $\mathcal L^0$, proves the growth rates are real, and bounds the fastest rate \cite{LinStrauss2008}. Their companion paper treats nonlinear stability and instability in the periodic setting \cite{LinStraussNonlinear2007}. Nguyen and Strauss obtained sharp criteria with specular particle reflection and conducting boundaries, first for radial and longitudinal symmetries and then in a solid torus \cite{NguyenStrauss2013,NguyenStrauss2014}. The second author extended this analysis to general axisymmetric conducting domains \cite{Zhang2019}. On unbounded domains Ben-Artzi and Holding gave sufficient instability criteria that allow nonmonotone equilibria \cite{BenArtziHolding2017}.

The use of velocity parity goes back to Antonov's work in stellar dynamics \cite{Antonov1961}. In plasma theory Kruskal and Oberman derived a guiding-center energy principle \cite{KruskalOberman1958}. Pfirsch used conservation of canonical momentum to obtain marginal equations for the cylindrical $z$- and $\theta$-pinches \cite[\S\S2.4--2.5]{Pfirsch1962}. Schindler studied perturbations constant along the equilibrium current, and Schindler, Pfirsch and Wobig developed a two-dimensional kinetic stability criterion for this symmetry \cite{Schindler1966,SchindlerPfirschWobig1973}. {\AA}gren gave a comparison theorem for two kinetic species and applied it to $z$-pinch pressure profiles \cite{Agren1986}. Symon, Seyler and Lewis formulated kinetic stability tests in terms of orbit invariants, with an explicit reduction by Lewis and Seyler when only one coordinate is nonignorable \cite{SymonSeylerLewis1982,LewisSeyler1982}. Guo used symmetry-dependent conservation laws to prove nonlinear stability for magnetic equilibria \cite{Guo1997,Guo1999}. Guo and Strauss also exhibited magnetically induced instability \cite{GuoStrauss2000}.

Morrison proposed a Poisson formulation of the nonlinear Vlasov--Maxwell equations \cite{Morrison1980}. Weinstein and Morrison identified an obstruction to the Jacobi identity in the original bracket \cite{Morrison1981}, and Marsden and Weinstein obtained a geometric bracket for divergence-free magnetic fields \cite{MarsdenWeinstein1982}. Holm, Marsden, Ratiu and Weinstein developed the energy--Casimir approach to nonlinear stability \cite{HolmMarsdenRatiuWeinstein1985}. Morrison and Pfirsch computed the free energy along dynamically accessible variations \cite{MorrisonPfirsch1989,MorrisonPfirsch1990}. Correa-Restrepo and Pfirsch found negative-energy perturbations of equilibria with an ignorable coordinate \cite{CorreaPfirsch1997}. For a Maxwellian background, Bell, Campos Pinto, Possanner and Sonnendr\"ucker used the linearized Hamiltonian structure to construct a structure-preserving discretization \cite{BellCamposPintoEtAl2026}.

Lin identified the density shift and constrained minimization for models of Lin and Strauss, using velocity parity and the separable Hamiltonian theory of Lin and Zeng \cite{Lin2022,LinZengSeparable2022}. Here the shift follows from the momentum conserved along a symmetry. In our nonrelativistic setting a Killing field $X_j$ generates a rigid motion preserving the equilibrium and the perturbations. For a species of mass $m_s$ and charge $q_s$, the energy $e_s$ and canonical momenta $p_{sj}$ give the equilibrium distribution in invariant potentials $\phi^0,A^0$ as
\begin{equation}\label{eq:intro-profile}
 F_s=\mu_s(e_s,p_{s1},\ldots,p_{sk}),\qquad
 e_s=\tfrac12 m_s|v|^2+q_s\phi^0,\qquad
 p_{sj}=(m_s v+q_s A^0)\cdot X_j .
\end{equation}
A change $A$ of vector potential changes $p_{sj}$ by $q_s A\cdot X_j$. Subtracting the resulting change of density exposes the velocity reflection that separates the Hamiltonian equations for the $z$- and $\theta$-pinches. In a screw-pinch, helical twist couples magnetic flux to the particle and electric constraints. At finite negative index, the general Hamiltonian theory of Lin and Zeng bounds the number of unstable eigenvalues \cite{LinZengIndex2022}. An additional symmetry restores separation in the zero helical harmonic.

The choice of momentum dependence distinguishes the equilibrium families. Komarov and Fadeev obtained the neutral multicomponent rigid rotor with a density proportional to $\sech^2(ar^2+b)$ \cite{KomarovFadeev1962}. Batchelor and Davidson treated kinetic $\theta$-pinch equilibria using angular momentum \cite{BatchelorDavidson1975}. Morse studied distributions with axial and angular drift and with quadratic dependence on canonical axial momentum \cite{Morse1969}, and Virtamo combined rotation and translation in locally Maxwellian screw-pinches \cite{Virtamo1980}. Cordonnier, Leoncini, Dif-Pradalier and Garbet used separate linear multipliers for axial and angular momentum \cite{Cordonnier2022}. Cordonnier, Lebouazda, Leoncini and Dif-Pradalier subsequently studied particle dynamics in these fields \cite{Cordonnier2026}. In an isothermal class, Crews obtained the Liouville $\theta$-pinch limit and related sheared flow to the velocity cumulants \cite{Crews2026}. A polynomial cumulant generating function has degree at most two by Marcinkiewicz's theorem \cite{Marcinkiewicz1939}. The Gaussian screw-pinch studied here depends quadratically on the single momentum $\ell+\kappa p_z$ \cite{HolmerZhang2026}.

Section~\ref{sec:framework} derives the constrained potential energy and the test for separation. Section~\ref{sec:z} computes the bare and accessible indices of the Bennett pinch and confines any $z$-independent growing mode of the variable-drift pinch to the radial harmonic. Section~\ref{sec:theta} proves axisymmetric spectral stability for the angular-momentum Gaussian under a bound on the slope at the axis and for the center-peaked rigid rotor. Section~\ref{sec:screw} locates the radial spectrum on the real and imaginary axes, proves that all sufficiently high helical harmonics have imaginary spectrum, and verifies spectral stability in the full helical symmetry class for two explicit Gaussian equilibria.

\subsection*{Disclosure of AI assistance}
In earlier notes dating back to 2022, we adapted Lin's separable Hamiltonian
formulation~\cite{Lin2022} and the Lin--Strauss stability
criteria~\cite{LinStrauss2007,LinStrauss2008} to $z$-independent
perturbations of the Bennett equilibrium~\cite{Bennett1934}.
Following our work on two-species equilibria~\cite{HolmerZhang2026},
we outlined the broader framework and derived the density correction
for translational, rotational, and helical symmetries from the Noether
identity in Lemma~\ref{lem:noether}.
We used ChatGPT (OpenAI) to help develop the proof of
Theorem~\ref{thm:separable}, which realizes the Lin--Zeng
self-adjoint square~\cite{LinZengSeparable2022} explicitly and
characterizes the entire nonzero spectrum 
without assuming a positive energy gap or finite negative index.
In the examples studied here, the theorem applies to $z$-independent
perturbations of the $z$-pinches, axisymmetric perturbations of the
$\theta$-pinches, and the zero helical harmonic of the screw-pinches.
ChatGPT also assisted with checking proofs and references and with
revising the manuscript. It developed and tested the independent
interval-arithmetic verifier used in Appendix~\ref{app:screw-certificate}.

\section{Split Hamiltonian structure for the Vlasov--Maxwell equations}\label{sec:framework}

Choose potentials independent of $z$. A particle conserves the axial
momentum $m_s v_z+q_s A_z$, so a change $a_z$ in the potential changes
$\mu_s(e_s,p_z)$ by $q_s\mu_{s,p_z}a_z$ at fixed particle energy.
Subtracting this known term from the distribution perturbation leaves
an equation driven only by the electric field. The magnetic field still
enters through the particle orbits and the perturbed current.

\subsection{Conserved momentum and the particle equation}
For finitely many species with masses $m_s>0$ and charges $q_s$, consider
\begin{align}
 \partial_t f_s+v\cdot\nabla_x f_s+\frac{q_s}{m_s}(E+v\times B)\cdot\nabla_v f_s&=0,\label{eq:VM}\\
 E_t=\curl B-j,\qquad B_t&=-\curl E,\label{eq:Maxwell}\\
 \diver E=\rho,\qquad \diver B&=0,\label{eq:Gauss}\\
 \rho=\sum_s q_s\int_{\R^3}f_s\dd v,\qquad
 j&=\sum_s q_s\int_{\R^3}v f_s\dd v.\nonumber
\end{align}
Let $(F_s,E^0,B^0)$ be stationary. From now on $(f_s,E,B)$ denotes a first-order perturbation. Set
\begin{equation}\label{eq:D}
 D_s=v\cdot\nabla_x+\frac{q_s}{m_s}(E^0+v\times B^0)\cdot\nabla_v.
\end{equation}
The linearized particle equation is
\begin{equation}\label{eq:linear-V}
 (\partial_t+D_s)f_s=-\frac{q_s}{m_s}(E+v\times B)\cdot\nabla_vF_s,
\end{equation}
and the perturbation fields obey \eqref{eq:Maxwell}--\eqref{eq:Gauss} with the linearized moments.

Let $X_1,\ldots,X_k$, $k=1$ or $2$, be commuting Killing fields that
preserve the domain and are independent on its regular set. If $\Phi_c$
is a symmetry, a particle density is invariant when
$f(\Phi_c(x),d\Phi_c(x)v)=f(x,v)$. For a rotation this fixes the
cylindrical velocity components as the spatial angle changes.
Potentials are invariant as one-forms, $\mathcal L_{X_j}A^\flat=0$.

\begin{assumption}[Noether data]\label{ass:noether}
The equilibrium has smooth invariant potentials $E^0=-\nabla\phi^0$,
$B^0=\curl A^0$ and solves the stationary Maxwell equations.
Its nonnegative profiles are $C^2$ functions of the energy and momenta
in \eqref{eq:intro-profile}. The fixed-flux perturbations admit invariant potentials
\begin{equation}\label{eq:potentials}
 E=-\nabla\phi-\partial_t A,\qquad B=\curl A.
\end{equation}
All velocity integrals below converge and integration by parts has zero boundary term.
\end{assumption}
\begin{lemma}[Noether identity]\label{lem:noether}
For a time-independent Killing field $X$ and invariant potentials,
\begin{equation}\label{eq:noether}
 X\cdot(E+v\times B)=-(\partial_t+v\cdot\nabla_x)(A\cdot X).
\end{equation}
The full particle quantity $m_s v\cdot X+q_s A\cdot X$ is conserved even for time-dependent invariant potentials. In particular,
\begin{equation}\label{eq:invariants}
 D_se_s=D_sp_{sj}=0,\qquad D_s\mu_{s,e}=D_s\mu_{s,p_j}=0,
 \qquad \nabla_vF_s=m_s\mu_{s,e}v+m_s\sum_j\mu_{s,p_j}X_j.
\end{equation}
\end{lemma}
\begin{proof}
Invariance gives $X_i\partial_i A_j=-A_i\partial_jX_i$ and $X\phi=0$.
Hence $X\cdot(v\times B)=-v\cdot\nabla(A\cdot X)$, which proves
\eqref{eq:noether}. Along the particle path the term
$m_sv_iv_j\partial_jX_i$ vanishes by the Killing equation. This gives
momentum conservation. Conservation of stationary particle energy and
the chain rule give \eqref{eq:invariants}.
\end{proof}

Let $P_\VV$ project onto $\spanop\{X_1,\ldots,X_k\}$ and $P_\HH=I-P_\VV$.
Thus $\VV$ and $\HH$ denote components along and perpendicular to the symmetry orbits. Reflect the perpendicular velocity and split the distribution into its even and odd parts,
\begin{equation}\label{eq:reflection}
 R_x=P_\VV-P_\HH,\qquad (\mathcal Rf)(x,v)=f(x,R_xv),\qquad P_\pm=(I\pm\mathcal R)/2.
\end{equation}
Write $v^\VV=P_\VV v$ and $v^\HH=P_\HH v$. Particle energy and the momenta $p_{sj}$ are even under this reflection, as are the profile derivatives.

Put $a=P_\VV A$. Invariance gives $a_t=-E^\VV$. Lin's Hamiltonian formulation of the periodic and axisymmetric models of Lin and Strauss uses a correction to the even part of the distribution \cite[\S\S3--4]{Lin2022}. For the momenta associated to $X_1,\ldots,X_k$, the correction takes the form
\begin{equation}\label{eq:shift}
 S_s[a]=\sum_j\mu_{s,p_j}(a\cdot X_j),\qquad
 g_s=P_+f_s-q_sS_s[a],\qquad h_s=P_-f_s.
\end{equation}
Recover the distribution from $f_s=g_s+h_s+q_sS_s[a]$.

\begin{proposition}[Common corrected equation]\label{prop:shift}
Under Assumption~\ref{ass:noether}, \eqref{eq:linear-V} is equivalent to
\begin{equation}\label{eq:corrected}
 (\partial_t+D_s)(g_s+h_s)=-q_s\mu_{s,e}\,v\cdot E.
\end{equation}
\end{proposition}
\begin{proof}
Substitution of \eqref{eq:invariants} in \eqref{eq:linear-V} gives the right side
\[
 -q_s\mu_{s,e}v\cdot E-q_s\sum_j\mu_{s,p_j}X_j\cdot(E+v\times B)
 =-q_s\mu_{s,e}v\cdot E+q_s(\partial_t+D_s)S_s[a],
\]
by Lemma~\ref{lem:noether}. Subtract the last term.
\end{proof}

On smooth invariant functions, the parity-preserving and parity-reversing parts of transport are
\begin{equation}\label{eq:Dsplit}
 D_{s,\pre}=\tfrac12(D_s+\mathcal RD_s\mathcal R),\qquad
 D_{s,\rev}=\tfrac12(D_s-\mathcal RD_s\mathcal R).
\end{equation}
Assume henceforth $w_s=-\mu_{s,e}>0$ on the occupied invariant phase space. The exact parity equations are
\begin{align}
 g_{s,t}&=-D_{s,\pre}g_s-D_{s,\rev}h_s+q_sw_s v^\VV\cdot E^\VV,\label{eq:g}\\
 h_{s,t}&=-D_{s,\rev}g_s-D_{s,\pre}h_s+q_sw_s v^\HH\cdot E^\HH.\label{eq:h}
\end{align}
The kinetic space and its pairing are
\begin{equation}\label{eq:kinspace}
 \cK=\bigoplus_s L^2(w_s^{-1}\dd x\dd v),\qquad
 \ip{f}{g}_\kin=\sum_s\iint f_s\cl{g_s}\,w_s^{-1}\dd x\dd v,
 \qquad \cK_\pm=P_\pm\cK.
\end{equation}
Here $\dd x$ is physical reduced volume, for example $r\dd r\dd\theta$ per unit axial length.

\subsection{Which symmetries separate the equations?}\label{sub:geometry}
For one nonvanishing Killing field write
\begin{equation}\label{eq:twist}
 X=\ell n,\qquad |n|=1,\qquad \tau=n\cdot\curl n.
\end{equation}
The Killing equation implies
\begin{equation}\label{eq:unit-generator}
 \nabla_n n=-\nabla\log\ell,\qquad
 \curl n=\tau n+\nabla\log\ell\times n,\qquad
 \nabla_u n=\tfrac\tau2 n\times u\quad(u\perp n).
\end{equation}
The last identity follows because the restriction of $\nabla n$ to $n^\perp$ is skew and its trace against rotation is $\tau$.

\begin{proposition}[Curl and transport under one symmetry]\label{prop:one-geometry}
On invariant fields, with the parallel component first,
\begin{equation}\label{eq:curl-one}
 \curl=\begin{pmatrix}\tau I&C^*\\ C&0\end{pmatrix},\qquad
 C(\alpha n)=\ell^{-1}\nabla(\ell\alpha)\times n.
\end{equation}
The star here denotes a formal adjoint in physical volume. On invariant distributions,
\begin{equation}\label{eq:Dpre-one}
 D_{s,\pre}=\left(\tau v_n+\frac{q_s}{m_s}B_n^0\right)
 (v^\HH\times n)\cdot\nabla_v.
\end{equation}
Consequently, for a charged species on an open velocity set, simultaneous off-diagonal curl and pure transport reversal hold precisely when $\tau=0$ and $B_n^0=0$.
\end{proposition}
\begin{proof}
For an invariant one-form $\beta$ perpendicular to $X$, Cartan's formula gives $\iota_Xd\beta=0$, since $\iota_X\beta=0$ and $\mathcal L_X\beta=0$. Its curl is therefore parallel to $X$. The product rule and \eqref{eq:unit-generator} give \eqref{eq:curl-one} for a field parallel to $X$. Integration by parts identifies the other block.

To compute transport, use $R=2n\otimes n-I$. With $v=v_n n+v^\HH$ and $\widetilde v=Rv$, differentiating $\mathcal R(v\cdot\nabla_x)\mathcal R$ produces the additional coefficient $(\nabla_{\widetilde v}R)\widetilde v$ in velocity. Tangent-lift invariance gives
\[
 n\cdot\nabla_x f=-\left(v_n\nabla_n n+\tfrac\tau2 n\times v^\HH
 -(v^\HH\cdot\nabla_n n)n\right)\cdot\nabla_vf.
\]
Taking the half-sum cancels the two terms containing $\nabla_n n$ and leaves $\tau v_n(v^\HH\times n)\cdot\nabla_v$. The equilibrium electric field is perpendicular to $X$. The parity-preserving magnetic term is precisely the rotation produced by $B_n^0n$, proving \eqref{eq:Dpre-one}. Its coefficient is affine in $v_n$, so vanishing on an open velocity set forces both coefficients to vanish.
\end{proof}

Since $j^0_\HH=0$, the stationary Maxwell equations give $C(B_n^0n)=0$. Consequently
\begin{equation}\label{eq:guide}
  \ell B_n^0=\text{constant}
\end{equation}
on each connected regular region. For translations and rotations $\tau=0$.
Translation permits a constant axial guide field. Rotation permits
$c/r$ in the azimuthal direction, and regularity at a solid axis forces
$c=0$. For $X=\partial_\theta+\kappa\partial_z$,
$\tau=2\kappa/(r^2+\kappa^2)$. Thus even a vanishing guide field leaves a
term that preserves particle parity.

The only pair of generators needed here consists of axial translation and
rotation. A doubly invariant field depends on $r$ alone. Its radial
component has zero curl and its tangential components have tangential curl,
\begin{equation}\label{eq:curl-two}
 \curl(a_\theta e_\theta+a_z e_z)
   =-a_z' e_\theta+r^{-1}(ra_\theta)'e_z=:C_\VV a .
\end{equation}
The invariant equilibrium potential gives $B^0_r=0$. In cylindrical
velocities, $v_r\partial_r$, the radial acceleration times $\partial_{v_r}$,
and the two tangential acceleration terms all reverse parity under
$v_r\mapsto-v_r$. Hence $D_{s,\pre}=0$. This observation gives a second
separation for the screw-pinch.

\subsection{The field energy and evolution}\label{sub:energy}\label{sub:formal}
Spatial fields use their usual $L^2$ pairing in physical reduced volume. The particle and field variables are coupled by
\begin{align}
 (\Gamma_\VV a)_s=q_sw_s v^\VV\cdot a,&\qquad
 (\Gamma_\HH e)_s=q_sw_s v^\HH\cdot e,\label{eq:Gamma}\\
 \Gamma_\VV^*g=\sum_s q_s\int v^\VV g_s\dd v,&\qquad
 \Gamma_\HH^*h=\sum_s q_s\int v^\HH h_s\dd v,\label{eq:Gamma-adj}\\
 Ka=\sum_s q_s^2\int v^\VV S_s[a]\dd v,&\qquad
 Qg=\sum_s q_s\int g_s\dd v.\label{eq:KQ}
\end{align}
The stars in \eqref{eq:Gamma-adj} are the adjoints for \eqref{eq:kinspace}. Parity gives
\begin{equation}\label{eq:moments}
 j=\Gamma_\VV^*g+\Gamma_\HH^*h+Ka,\qquad
 \rho=Q(g+\Gamma_\VV a).
\end{equation}

\begin{lemma}[Current coefficient]\label{lem:K}
With the velocity integration by parts in Assumption~\ref{ass:noether},
\begin{equation}\label{eq:Kidentity}
 K=\Gamma_\VV^*\Gamma_\VV-d(x)I_\VV,\qquad
 d(x)=\sum_s\frac{q_s^2}{m_s}\int F_s(x,v)\dd v.
\end{equation}
In particular $K(x)$ is symmetric. Every invariant field tangent to the symmetry orbits is divergence-free, so Gauss's law is
\begin{equation}\label{eq:Gauss-shift}
 \diver E^\HH=Q(g+\Gamma_\VV a).
\end{equation}
\end{lemma}
\begin{proof}
For $a$ tangent to the symmetry orbits, \eqref{eq:invariants} gives
$S_s[a]=w_s v\cdot a+m_s^{-1}a\cdot\nabla_vF_s$.
Integrating its product with $v\cdot c$ gives the second term
$-m_s^{-1}(a\cdot c)\int F_s\dd v$, proving \eqref{eq:Kidentity}. Integrating without $v\cdot c$ proves the charge identity in \eqref{eq:moments}. Finally write an invariant field tangent to the symmetry orbits as $Z=\sum_j a_jX_j$. Commutation and independence of the generators give $X_i a_j=0$, and each $X_j$ is divergence-free. Hence $\diver Z=0$.
\end{proof}

\begin{proposition}[Conserved quadratic energy]\label{prop:energy}
For finite-energy classical solutions with vanishing boundary forms,
\begin{equation}\label{eq:energy}
 2\mathcal H_2=\norm g_\kin^2+\norm h_\kin^2+\norm E^2+\norm B^2-\ip{Ka}{a}
\end{equation}
is constant.
\end{proposition}
\begin{proof}
Pair \eqref{eq:corrected} with $g+h$ and use skew transport. The kinetic
energy derivative is $\operatorname{Re}\ip{E}{j-Ka}$, and Maxwell's
equations contribute $-\operatorname{Re}\ip{E}{j}$. Since
$a_t=-E^\VV$, the derivative of $-\frac12\ip{Ka}{a}$ cancels the difference.
\end{proof}
For one symmetry take $b=B^\VV$, the total component along the generator.
Then
\begin{equation}\label{eq:mag-coordinates}
 B=Ca+b,\qquad b-\tau a=C^*A^\HH,\qquad
 M=C^*C-K.
\end{equation}
For two symmetries $B=C_\VV a$, and $M=C_\VV^*C_\VV-K$.
In either case the potential contributes $m[a]=\|C_ea\|^2-\ip{Ka}{a}$
to the energy, where $C_e$ is the realization of the appropriate curl
specified below.

For one generator group the variables as
\begin{equation}\label{eq:uv}
 u=(g,E^\HH,a)^T,\qquad \upsilon=(h,E^\VV,b)^T,\qquad
 L_u=\operatorname{diag}(I,I,M).
\end{equation}
The six equations are
\begin{align}
 g_t&=-D_\pre g-D_\rev h+\Gamma_\VV E^\VV,&
 E_t^\HH&=-\Gamma_\HH^*h+Cb,& a_t&=-E^\VV,\label{eq:system-one-a}\\
 h_t&=-D_\rev g-D_\pre h+\Gamma_\HH E^\HH,&
 E_t^\VV&=-\Gamma_\VV^*g+Ma+\tau b,&
 b_t&=-C^*E^\HH-\tau E^\VV.\label{eq:system-one-b}
\end{align}
Thus
\begin{equation}\label{eq:general-matrix}
 \partial_t\binom u\upsilon
 =\begin{pmatrix}G_+&\cB\\-\cB^*&G_-\end{pmatrix}
  \begin{pmatrix}L_u&0\\0&I\end{pmatrix}\binom u\upsilon,
\end{equation}
where in these physical coordinates the diagonal blocks are
\begin{equation}\label{eq:Gformal}
 G_+u=(-D_\pre g,0,0),\qquad
 G_-\upsilon=(-D_\pre h,\tau b,-\tau E^\VV).
\end{equation}
The off-diagonal blocks are given by \eqref{eq:system-one-a}--\eqref{eq:system-one-b}.
We call the system \emph{separable} when $G_+=G_-=0$. In the radial
two-generator case we omit $b$ and obtain
\begin{equation}\label{eq:two-matrix}
 u_t=\cB\upsilon,\qquad \upsilon_t=-\cB^* L_u u,\qquad
 \cB(h,e)=(-D h+\Gamma_\VV e,-\Gamma_\HH^* h,-e).
\end{equation}

\subsection{The closed evolution and its constraints}\label{sub:closed}\label{sub:accessibility}

Write $H_\VV,H_\HH$ for the spatial $L^2$ spaces of invariant fields in the
two directions. We specify how transport and curl act on the energy space,
since the location of the full spectrum depends on their domains.
\begin{assumption}[Flow and couplings]\label{ass:analytic}
For each species, the complete equilibrium characteristic flow, including the chosen wall reflection, exists almost everywhere for all positive and negative times. It preserves phase volume, the occupied region, and $w_s$. Writing $\Phi_s^t$ for this flow, $U_s(t)k=k\circ\Phi_s^{-t}$ defines a strongly continuous unitary group on the weighted kinetic space and preserves the retained symmetry class. Its generator is $-D_s$. At a wall we require the reflection to preserve particle energy and the retained canonical momenta. Specular reflection does so on every wall tangent to the symmetry generators.

The maps $\Gamma_\VV,\Gamma_\HH,Q$ are bounded on the indicated $L^2$ spaces, $K$ is bounded self-adjoint, and $\tau$ is a bounded real multiplier when present. To recover the original density completion, also assume
\begin{equation}\label{eq:Sigma-bound}
 \Sigma:a\longmapsto(q_sS_s[a])_s
 \quad\hbox{is bounded from the potential graph space into }\cK_+.
\end{equation}
\end{assumption}

A sufficient set of moment bounds for the current and charge maps is
\begin{equation}\label{eq:moment-bounds}
 \esssup_x\sum_s q_s^2\int(1+|v|^2)w_s\dd v<\infty,
 \qquad
 \esssup_x\sum_s\frac{q_s^2}{m_s}\int F_s\dd v<\infty.
\end{equation}
The shift bound \eqref{eq:Sigma-bound} is separate. If weights vanish outside an occupied region, its invariant boundary and the allowed perturbations are part of the flow hypothesis.

For one generator, $C_e$ and $C_m$ are closed densely defined realizations of $C$ from $H_\VV$ into $H_\HH$. At a conducting wall in the fixed-potential sector we take $C_e=C_{\min}$ and $C_m=C_{\max}$. On their respective domains, $a$ and $E^\VV$ have zero wall trace, $C_m^*$ imposes the vanishing tangential trace of $E^\HH$, and $b$ has the maximal magnetic domain. For two generators the curl in \eqref{eq:curl-two} has two tangential components, so $C_e:D(C_e)\subset H_\VV\to H_\VV$ realizes $C_\VV$, with the minimal domain at a conducting wall. The variables $b,C_m$ then drop out. All closures extend smooth
Cartesian invariant fields across the axis. Periodic coordinates are
identified before completion. Put
\begin{equation}\label{eq:Lambda}
 \cV=D(C_e),\qquad
 \norm a_\cV^2=\norm a^2+\norm{C_ea}^2,\qquad
 \Lambda=(I+C_e^*C_e)^{1/2},\qquad \widehat a=\Lambda a.
\end{equation}
The graph coordinate is unitary from $\cV$ onto $H_\VV$. On the whole cylinder the condition $a\in L^2$ is part of the state space. A perturbation with square-integrable magnetic field can have a symmetry potential outside this space. All full-spectrum assertions below refer to the graph realization \eqref{eq:Lambda}. The potential form and its normalized bounded energy operator are
\begin{align}
 m[a,c]&=\ip{C_ea}{C_ec}-\ip{Ka}{c},&
 M&=C_e^*C_e-K,\label{eq:potential-form}\\
 L_a&=I-\Lambda^{-1}(I+K)\Lambda^{-1},&
 m[a,c]&=\ip{L_a\widehat a}{\widehat c}.\label{eq:La}
\end{align}

Use the normalized state and energy
\begin{equation}\label{eq:normalized-space}
 z=(g,E^\HH,\widehat a;h,E^\VV,b)\in\cE,\qquad
 \cL=\operatorname{diag}(I,I,L_a,I,I,I),
\end{equation}
where $\cE$ is the product of the kinetic and spatial $L^2$ spaces just
defined. Appendix~\ref{app:framework} gives the skew operator $\cJ$, its
domain, and the proof of the following theorem.
\begin{theorem}[Closed Hamiltonian evolution]\label{thm:closed}
Under Assumptions~\ref{ass:noether} and~\ref{ass:analytic} and the chosen field realizations, $\cJ^*=-\cJ$. The natural product
\begin{equation}\label{eq:Aclosed}
 \cA=\cJ\cL,\qquad D(\cA)=\{z:\cL z\in D(\cJ)\}=D(\cJ)
\end{equation}
generates a strongly continuous group on $\cE$ and preserves the bounded energy form
$q(z,w)=\ip{\cL z}{w}$. It realizes \eqref{eq:system-one-a}--\eqref{eq:system-one-b}, or \eqref{eq:two-matrix}, with the domains in Appendix~\ref{app:framework}. Under \eqref{eq:Sigma-bound}, the inverse density shift identifies it with the augmented weak linearized Vlasov--Maxwell system.
\end{theorem}

A split is separable at the level of closed operators if $\tau=0$ and
the reflected characteristic flow runs backward,
\begin{equation}\label{eq:strong-reversal}
 \mathcal R U(t)\mathcal R=U(-t).
\end{equation}
It implies invariance of $D(D)$ under $\mathcal R$ and gives closed restrictions $D_-:\cK_-\to\cK_+$, $D_+:\cK_+\to\cK_-$ with $D_-^*=-D_+$. The closed coupling is
\begin{equation}\label{eq:Bclosed}
 \cB(h,e,b)=(-D_-h+\Gamma_\VV e,-\Gamma_\HH^*h+C_mb,-\Lambda e),
 \quad D(\cB)=D(D_-)\oplus D(\Lambda)\oplus D(C_m),
\end{equation}
with
\begin{equation}\label{eq:Bclosedadj}
 \cB^*(g,E,\widehat a)=(D_+g-\Gamma_\HH E,\Gamma_\VV^*g-\Lambda\widehat a,C_m^*E).
\end{equation}
The adjoint domain is $D(D_+)\oplus D(C_m^*)\oplus D(\Lambda)$. After subtracting the bounded current maps, the derivatives act on separate slots. This gives the closed coupling used in the separable theory \cite{Lin2022,LinZengSeparable2022}.

Let $\cT$ denote the invariant scalar $H^1$ test functions, periodic in periodic directions and of zero trace on conducting walls. The field realizations above give $\nabla\chi\in D(C_m^*)$ and $C_m^*\nabla\chi=0$ for $\chi\in\cT$ when $b$ is present. For other realizations we require these two properties to impose the same Gauss tests. The particle and Maxwell equations conserve both Gauss's law and magnetic
flux. More generally $\ip{z(t)}{n}$ is constant for every $n\in\Ker\cJ$.
We call
\begin{equation}\label{eq:X}
 \cX=\cl{\Ran\cJ}=(\Ker\cJ)^\perp
\end{equation}
the accessible space for the linearized evolution and write $\cL_\cX=P_\cX\cL|_\cX$ for its energy operator. Equivalently, its elements have zero pairing with each of these conserved linear functionals. Within the same augmented space, $\cY$ imposes Gauss's law and fixed magnetic flux. Appendix~\ref{app:framework} shows $\cX\subset\cY$ and
\begin{equation}\label{eq:nonzero-spectrum}
 \sigma(\cA)\setminus\{0\}=\sigma(\cA_\cY)\setminus\{0\}
 =\sigma(\cA_\cX)\setminus\{0\}.
\end{equation}
Thus the generalized eigenspace at every $\lambda\ne0$ is contained in $\cX$. Appendix~\ref{app:framework} also shows that a physical normal mode at such a $\lambda$ with finite kinetic and Maxwell norms, in the corresponding characteristic and Maxwell domains, lifts to this graph realization.

\subsection{Reducing the energy to the potential}\label{sub:Schur}
Lin computed a constrained minimum at fixed potential for the separable Vlasov--Maxwell models \cite[\S\S3--4]{Lin2022}. Here the same minimization includes the magnetic flux when the symmetry has twist. Let $\cF=\cK\oplus H_\HH\oplus H_\VV$ and let $P_\cZ$ project
orthogonally onto
\begin{equation}\label{eq:Zall}
 \cZ=\left\{(k,e,b):
 \begin{array}{l}
 k\in D(D),\ e\in D(C_m^*),\ b\in D(C_m),\\
 Dk=\Gamma_\HH e,\quad C_mb=\Gamma_\HH^*P_-k,\quad C_m^*e=0
 \end{array}\right\}.
\end{equation}
For two generators omit $b,C_m,C_m^*$, and replace the middle equation by $\Gamma_\HH^*P_-k=0$. There is then no derivative domain condition on $e$. This is a closed joint kernel because the three unbounded derivatives have separate input slots and the couplings are bounded. Set
\begin{equation}\label{eq:Gammaall}
 \Gamma_{\rm all}a=(\Gamma_\VV a,0,-\tau a)\in\cF.
\end{equation}
\begin{theorem}[Accessible spatial energy]\label{thm:accessibility}
The effective potential form and spatial operator are
\begin{align}
 m_\acc[a]&=\norm{C_ea}^2-\ip{Ka}{a}
              +\norm{P_\cZ\Gamma_{\rm all}a}_{\cF}^2,\qquad D(m_\acc)=\cV,\label{eq:macc}\\
 M_\acc&=C_e^*C_e-K+\Gamma_{\rm all}^*P_\cZ\Gamma_{\rm all}.
 \label{eq:Macc}
\end{align}
The operator is self-adjoint on $D(C_e^*C_e)$ and bounded below. Its negative index satisfies
\begin{equation}\label{eq:accessible-index}
 N_\acc:=\nm(\cL_\cX)=\nm(m_\acc)=\nm(M_\acc)
 \le\nm(m).
\end{equation}
The indices may be infinite. No inverse of $D$ or closed transport range is assumed.
\end{theorem}
\begin{proof}
Solving \eqref{eq:Jclosed} for a Poisson-kernel vector gives $E^\VV=0$, $(k,e,b)\in\cZ$, and
\[
 g=P_+k,\quad h=P_-k,\quad
 \widehat a=\Lambda^{-1}(\Gamma_\VV^*P_+k-\tau b).
\]
The last expression lies in $D(\Lambda)$, verifying the remaining domain condition. Orthogonality to this graph says, in physical coordinates $p=(g+h,E^\HH,b)$,
\[
 p+\Gamma_{\rm all}a\in\cZ^\perp.
\]
At fixed $a$ the unique minimizing $p$ is $-P_\cZ\Gamma_{\rm all}a$, and the minimizing parallel electric field is zero. Write $w=p+\Gamma_{\rm all}a$. Orthogonal completion of $\norm p^2+\norm{E^\VV}^2+m[a]$ gives
\begin{equation}\label{eq:completed-square}
 q(z,z)=\norm{w-P_{\cZ^\perp}\Gamma_{\rm all}a}^2
       +\norm{E^\VV}^2+m_\acc[a].
\end{equation}
This is a bounded invertible congruence between $\cX$ and
$\cZ^\perp\oplus\cV\oplus H_\VV$, with its graph norm. The positive squares prove the index equality. The correction in \eqref{eq:Macc} is bounded and nonnegative, so the form representation theorem gives the asserted self-adjoint operator and the remaining inequality.
\end{proof}

Here $n^-$ denotes the dimension of the negative spectral subspace,
equivalently the supremum of the dimensions of strictly negative subspaces.

Suppose \eqref{eq:strong-reversal} holds and $\tau=0$. Define
\begin{equation}\label{eq:Zsep}
 \cZ_+=\{(g,e):g\in D(D_+),\ e\in D(C_m^*),\ D_+g=\Gamma_\HH e,
                     \ C_m^*e=0\}.
\end{equation}
The full joint kernel splits into even/electric and odd/magnetic parts. Since $\Gamma_{\rm all}a$ is even and has no magnetic component,
\begin{equation}\label{eq:maccsep}
 m_\acc[a]=m[a]+\norm{P_{\cZ_+}(\Gamma_\VV a,0)}^2.
\end{equation}
In particular a projection onto guessed particle invariants alone generally misses an electric correction.

We can now recover the scalar operator of Lin and Strauss. Let $\Pi=P_{\Ker D_+}$ and put
\begin{equation}\label{eq:scalarform}
 T_\phi\chi=(I-\Pi)Q^*\chi,\qquad
 a_\phi(\chi,\psi)=\ip{\nabla\chi}{\nabla\psi}+\ip{T_\phi\chi}{T_\phi\psi}.
\end{equation}
Complete $\cT$ modulo this form's nullspace. The map
$\chi\mapsto(T_\phi\chi,\nabla\chi)$ is an isometry onto a closed subspace $\cZ_\nabla\subset\cZ_+$, orthogonal to $\cI_0=(\Ker D_+)\oplus\{0\}$. Retain
$\cZ_{\rm rem}=\cZ_+\ominus(\cI_0\oplus\cZ_\nabla)$. Let $\chi_a$ be the Riesz solution
\begin{equation}\label{eq:scalar-solve}
 a_\phi(\chi_a,\psi)=\ip{\Gamma_\VV a}{(I-\Pi)Q^*\psi}.
\end{equation}
Orthogonal projection onto the three summands gives
\begin{equation}\label{eq:Schursep}
 m_\acc[a]=m[a]+\norm{\Pi\Gamma_\VV a}^2
       +a_\phi(\chi_a,\chi_a)+\norm{P_{\cZ_{\rm rem}}(\Gamma_\VV a,0)}^2.
\end{equation}
Equivalently, by \eqref{eq:Kidentity}, the first two terms are
$\norm{C_ea}^2+\ip{d a}{a}-\norm{(I-\Pi)\Gamma_\VV a}^2$.

Writing the field amplitudes as scalars gives the Lin--Strauss operators \cite{LinStrauss2007,LinStrauss2008} for one generator,
\begin{align}
 \mathcal A_1^0&=-\Delta+Q(I-\Pi)Q^*,&
 \mathcal B^0&=Q(I-\Pi)\Gamma_\VV,\label{eq:LS-one}\\
 \mathcal A_2^0&=M+\Gamma_\VV^*\Pi\Gamma_\VV,&
 \mathcal L^0&=(\mathcal B^0)^*(\mathcal A_1^0)^{-1}
                   \mathcal B^0+\mathcal A_2^0 .\label{eq:LS-schur}
\end{align}
In particular $\mathcal A_1^0$ represents $a_\phi$, and
$\chi_a=(\mathcal A_1^0)^{-1}\mathcal B^0a$ means the weak solve
\eqref{eq:scalar-solve} on the completed potential space modulo its
nullspace. This convention also applies on an unbounded cross section.
If $\cZ_{\rm rem}=0$, then $M_\acc=\mathcal L^0$. In general the full
quadratic form is $\ip{\mathcal L^0a}{a}$ plus the final squared norm in
\eqref{eq:Schursep}. With two generators $a$ has two components and the
same elimination gives a matrix operator on them. On an even kinetic input $f$ the projection $\Pi f$ is the strong limit
of $T^{-1}\int_0^T U(s)f\,ds$ as $T\to\infty$, by the mean ergodic
theorem. This gives it directly in terms of equilibrium particle orbits.

\subsection{Location of the spectrum}\label{sub:spectral}
The separable theory of Lin and Zeng \cite{LinZengSeparable2022} leads to a self-adjoint second-order equation. In our graph coordinates its operator differs from $\cB^*\cB$ by a bounded spatial multiplier. This formula gives a resolvent argument even when the positive spectrum of the energy accumulates at zero or the negative index is infinite.
\begin{theorem}[Self-adjoint square and separable count]\label{thm:separable}
For the closed separable realization, let $L_U=\operatorname{diag}(I,I,L_a)$, let $\Pi_E$ extract $E^\VV$ from the second block, and set
\begin{equation}\label{eq:S}
 S=\cB^*\cB-\Pi_E^*(I+K)\Pi_E,\qquad D(S)=D(\cB^*\cB).
\end{equation}
Then $S$ is self-adjoint, equals the natural product $\cB^*L_U\cB$, and
\begin{equation}\label{eq:spectral-correspondence}
 \lambda\in\sigma(\cA),\quad\lambda\ne0
 \quad\Longleftrightarrow\quad -\lambda^2\in\sigma(S).
\end{equation}
Consequently $\sigma(\cA)\subset\R\cup i\R$ without compactness or finite index. Moreover
\begin{equation}\label{eq:sharp-count}
 \nm(S)=N_\acc=\nm(M_\acc).
\end{equation}
Thus the full spectrum is imaginary if and only if $m_\acc\ge0$. If $N_\acc<\infty$, the positive real eigenvalues have total algebraic multiplicity $N_\acc$ and there is no off-axis spectrum. All nonzero eigenvalues have no nontrivial Jordan chains. A nonzero imaginary eigenvector $(u,\upsilon)$ has energy $q((u,\upsilon),(u,\upsilon))=2\norm\upsilon^2>0$.
\end{theorem}
\begin{proof}
The potential output of $\cB\upsilon$ is $-\Lambda\Pi_E\upsilon$. Therefore on $D(\cB)$,
\begin{equation}\label{eq:square-form}
 \ip{L_U\cB\upsilon}{\cB\omega}
 =\ip{\cB\upsilon}{\cB\omega}
 -\ip{(I+K)\Pi_E\upsilon}{\Pi_E\omega}.
\end{equation}
This is a bounded perturbation of a closed nonnegative form, with operator \eqref{eq:S}. Also $(I-L_U)\cB\upsilon\in D(\cB^*)$ and its image under $\cB^*$ is $\Pi_E^*(I+K)\Pi_E\upsilon$. This verifies the natural product domain.

For $\lambda\ne0$, the resolvent system is
\[
 u=\lambda^{-1}(f+\cB\upsilon),\qquad
 (S+\lambda^2)\upsilon=\lambda g-\cB^*L_Uf,
\]
Here $(f,g)$ belongs to the two energy blocks. The second equation is a form equation on $D(\cB)$, where $\cB^*L_U f$ acts on $\omega$ by $\ip{L_U f}{\cB\omega}$. Choose $c$ so that $S+cI\ge I$. Relative to $(S+cI)^{1/2}$, its form operator has multiplier $(s+\lambda^2)/(s+c)$, with a bounded inverse exactly when $-\lambda^2\notin\sigma(S)$. This produces $\upsilon\in D(\cB)$, $u$ in the first energy block, and $L_Uu\in D(\cB^*)$, with a bounded resolvent. Conversely, a normalized approximate eigenvector $\upsilon_j\in D(S)$ at $-\lambda^2$ gives $(\lambda^{-1}\cB\upsilon_j,\upsilon_j)$ with resolvent residual tending to zero and norm at least one. Normalizing these pairs proves \eqref{eq:spectral-correspondence} for the entire nonzero spectrum.

The negative index of \eqref{eq:square-form} is the index of $L_U$ on $\Ran\cB$, and equals its index on $\cl{\Ran\cB}$ by approximation of any finite strictly negative Gram matrix. The accessible space is $\cl{\Ran\cB}\oplus\cl{\Ran\cB^*}$. The second energy block is positive. Theorem~\ref{thm:accessibility} now proves \eqref{eq:sharp-count}. At finite index the negative spectrum of $S$ consists of finitely many eigenvalues. The map $\upsilon\mapsto(\lambda^{-1}\cB\upsilon,\upsilon)$ identifies their eigenspaces with those of $\cA$ at $\lambda>0$.

A Jordan chain at nonzero $\lambda$ would give
$(S+\lambda^2)\upsilon_1=-2\lambda\upsilon_0$ and
$(S+\lambda^2)\upsilon_0=0$, contradicting self-adjointness after pairing with $\upsilon_0$. At $\lambda=i\omega$, \eqref{eq:square-form} gives $\ip{L_Uu}{u}=\norm\upsilon^2$, proving the energy statement. The physical constraint and accessible restrictions have identical nonzero spectral data by \eqref{eq:nonzero-spectrum}.
\end{proof}

\paragraph{Recovering growth rates.}
Lin and Strauss bounded the maximal growth rate by the spatial energy operator in their sharp criterion \cite{LinStrauss2008}. There is also a spatial equation for each positive rate in the present separable system.
It uses a transport solve at each trial value and thus provides a way to
compute more than the count. Write $e=E^\VV$, $y=(h,b)$, and set
\[
 T y=(-D_-h,-\Gamma_\HH^*h+C_m b),\qquad Ge=(\Gamma_\VV e,0).
\]
Omit $b$ in the two-generator case. The operator $T$ is closed on the
product derivative domain because its only mixed term is bounded.
By \eqref{eq:square-form}, the form of $S+tI$ is
\[
  \|Ty+Ge\|^2+t\|y\|^2+m[e]+t\|e\|^2.
\]
For $t>0$ minimize over $y$ by the normal equation. Its unique solution
is $y_t=-T^*(TT^*+tI)^{-1}Ge$. Substitution gives the spatial operator
\begin{equation}\label{eq:growth-pencil}
  \mathcal F(t)=M+tI+G^*t(TT^*+tI)^{-1}G.
\end{equation}
The same variation in $e$ shows that $\Ker\mathcal F(t)$ corresponds
bijectively to $\Ker(S+tI)$ by $e\mapsto(y_t,e)$ in this order of
coordinates. Thus $t=\lambda^2$ for a growing mode precisely when
$\mathcal F(t)$ has a zero eigenvalue. By the spectral theorem, for
$t_2>t_1>0$ one has $\mathcal F(t_2)-\mathcal F(t_1)\ge(t_2-t_1)I$.
As $t\downarrow0$, the multiplier $t(TT^*+tI)^{-1}$ tends strongly to
$P_{\Ker T^*}=P_{\cZ_+}$. Hence $\mathcal F(0+)=M_\acc$ in the sense of
forms. In particular any squared growth rate satisfies
$0<t\le-\inf\sigma(M_\acc)$. With finite negative index, one can search
this interval for zero crossings of a self-adjoint spatial problem.
The imaginary continuous spectrum is a separate question governed by
the nonnegative spectrum of $S$ in \eqref{eq:spectral-correspondence}.

\begin{theorem}[General accessible spectral bounds]\label{thm:general-spectral}
On the complexification of the real perturbation space the following bounds also apply to nonzero helical twist.
\begin{enumerate}[label=(\roman*),leftmargin=2em]
\item If $m_\acc\ge0$, then $\sigma(\cA)\subset i\R$. The group on $\cX$ has at most linear growth; the full group has at most quadratic growth. If already $m\ge0$, the full group has at most linear growth.
\item If $N_\acc<\infty$, all nonimaginary spectrum consists of finitely many isolated eigenvalues of finite algebraic multiplicity. They have the Hamiltonian quartet symmetry, with equal multiplicities. Writing $k_r$ for the positive-real algebraic count and $k_c$ for the first-quadrant algebraic count,
\begin{equation}\label{eq:unstable-bound}
 k_r+2k_c\le N_\acc.
\end{equation}
In particular $N_\acc\le1$ implies that the full spectrum lies on the axes.
\end{enumerate}
\end{theorem}
See Appendix~\ref{app:framework} for the proof.
An eigenvalue strictly off both axes requires at least two negative directions. Appendix~\ref{app:framework} also records the
stronger index identity used for individual low Fourier pairs in Section~\ref{sec:screw}.

\section{The \texorpdfstring{$z$}{z}-pinch}\label{sec:z}

In a $z$-pinch the particles carry current along the axis and produce a magnetic field circling it. We consider the classical Bennett equilibrium \cite{Bennett1934} and the variable-drift equilibrium studied in \cite[Theorems~1.1--1.2]{HolmerZhang2026}. The perturbations are independent of $z$ and may vary in both transverse coordinates. Spatial energies are per unit axial length. Schindler, Pfirsch and Wobig \cite{SchindlerPfirschWobig1973} derived kinetic stability criteria and a growth-rate expression for perturbations constant along the equilibrium current. {\AA}gren \cite{Agren1986} treated two kinetic species with profiles depending on energy and axial canonical momentum, including applications to $z$-pinch pressure profiles. Here we evaluate the accessible energy of Section~\ref{sec:framework} in this geometry.

\subsection{The operator for the axial potential}
For use in all three geometries, the cylindrical field relations and canonical momenta are
\begin{align}
 E^0&=-\phi^{0\prime}\eR,\qquad
 A^0=A_\theta^0\eTh+A_z^0\eZ,\qquad
 B_\theta^0=-A_z^{0\prime},\quad B_z^0=(rA_\theta^0)'/r,\label{eq:radial-background}\\
 \ell_s&=r(m_sv_\theta+q_s A_\theta^0),\qquad
 p_{sz}=m_sv_z+q_sA_z^0,\label{eq:cyl-momenta}\\
 -(r\phi^{0\prime})'/r&=\rho^0,\qquad
 -(B_z^0)'=j_\theta^0,\qquad (rB_\theta^0)'/r=j_z^0.\label{eq:radial-Maxwell}
\end{align}
Set $A_\theta^0=B_z^0=0$ and take $F_s=\mu_s(e_s,p_{sz})$, $w_s=-\mu_{s,e}>0$, where $e_s=m_s|v|^2/2+q_s\phi^0$. Since axial translation is the retained symmetry, reflection sends $(v_r,v_\theta,v_z)$ to $(-v_r,-v_\theta,v_z)$. It reverses the equilibrium transport on $z$-independent distributions. Indeed the transverse velocity changes sign, the transverse magnetic acceleration depends on $v_z$, and the axial acceleration is proportional to $v_r$. In Cartesian coordinates this proves the flow reversal through the axis as well.

Write $a=A_z$ for the perturbation, $g_s=f_{s,\ev}-q_s\mu_{s,p_z}a$, and $h_s=f_{s,\od}$. The coefficients entering the energy in Section~\ref{sec:framework} are
\begin{equation}\label{eq:z-KQ}
 (\Gamma_z a)_s=q_sw_sv_z a,\qquad
 K_z=\sum_s q_s^2\int v_z\mu_{s,p_z}\dd v
      =\partial_{A_z}j_z^0 .
\end{equation}
The derivative holds $r$ and $\phi^0$ fixed. Since $\curl(a\eZ)=(r^{-1}\partial_\theta a)\eR-(\partial_ra)\eTh$, the bare magnetic operator and the accessible form are
\begin{align}
 M_z&=-\Delta_\perp-K_z,\quad
 m_z[a]=\int_{\R^2}(|\nabla_\perp a|^2-K_z|a|^2)\dd x_\perp,\label{eq:z-M}\\
 \ip{M_{z,\acc} a}{a}_{\rm form}
   &=m_z[a]+\|P_{\cZ_{z,+}}(\Gamma_z a,0)\|^2 .\label{eq:z-macc}
\end{align}
Here $\cZ_{z,+}$ consists of pairs $(k,E_\perp)$ in the domains from Section~\ref{sec:framework} satisfying $D_+k=(q_s w_s v_\perp\cdot E_\perp)_s$ and $\curl E_\perp=0$. Thus $M_{z,\acc}$ includes the particle invariants and Gauss's law. On the whole plane its form domain is $H^1(\R^2)$. A circular conducting wall gives the same formulas on the disk with $a=0$ at the wall, specular reflection, and zero perturbation of the total axial magnetic flux.

\begin{theorem}[$z$-independent spectrum]\label{thm:z-application}
For a smooth stationary $z$-pinch of the above form, assume the flow, moment and density-shift hypotheses of Section~\ref{sec:framework}. The spectrum of the linearized evolution on the phase space with potential norm \eqref{eq:Lambda} is contained in $\R\cup i\R$. It is contained in $i\R$ precisely when $M_{z,\acc}\geq 0$. If $n^-(M_{z,\acc})<\infty$, this index counts the positive real eigenvalues with algebraic multiplicity. These conclusions also hold on the constrained and accessible spaces.
\end{theorem}
\begin{proof}
Translation has zero twist and $B_z^0=0$. Reflection therefore gives the off-diagonal kinetic blocks, and the axial curl and current above give the field blocks of the separable system. Theorems~\ref{thm:accessibility} and~\ref{thm:separable} apply with \eqref{eq:z-macc}.
\end{proof}

For the following applications, $j\in\mathbb Z$ denotes angular Fourier number. Expanding $a=\sum_j a_j(r)e^{ij\theta}/\sqrt{2\pi}$ defines
\begin{equation}\label{eq:z-Hj}
 M_z(a_je^{ij\theta})=(H_ja_j)e^{ij\theta},\qquad
 H_j=-\partial_r^2-r^{-1}\partial_r+j^2r^{-2}-K_z(r).
\end{equation}
Here $H_j$ is the self-adjoint operator on $L^2(r\dd r)$ obtained from the form of $M_z$ on the $j$th harmonic. The form domain consists of the radial coefficients of $H^1$ functions in this harmonic. These scalar operators include the equilibrium current derivative and describe the potential form before the positive correction in \eqref{eq:z-macc}.

For the four nontrivial families treated below in Sections~\ref{sec:z} and~\ref{sec:theta}, the full accessible operator equals the Lin--Strauss scalar operator \eqref{eq:LS-schur}. We record the reason, since the radial domain is unbounded. In each example $W(r)=\sum_s q_s^2\int w_s\dd v=O((1+r)^{-4})$, as follows from the densities below. For the variable-drift $z$-pinch, $(rU')'>0$ on the positive branch gives $U(r)\ge c_1\log r$ for large $r$ and some $c_1>0$. A curl-free transverse field on $\R^2$ is a gradient $\nabla\chi$. On $\R^2\times\T$ its possible axial period is constant in $r$ and vanishes by square integrability, so the same conclusion holds. In either geometry $\chi\in H^1_{\rm loc}$ can be chosen with $\int W|\chi|^2<\infty$. Indeed, on successive transverse annuli $2^j<r<2^{j+1}$, Poincar\'e's inequality on pairs of adjacent annuli bounds the difference of neighboring means by $C\|\nabla\chi\|_2$. The mean of $|\chi|^2$ on the $j$th annulus is thus $O((1+j)^2)$, sufficient for the weighted integral.

Truncate this potential at heights $\pm N$, multiply it by logarithmic cutoffs, then smooth. The approximants converge in $L^2(W\dd x)$ and in gradient norm. Indeed, for fixed $N$, the squared gradient cost of the cutoff between $R$ and $R^2$ is $O(N^2/\log R)$. Hence $Q^*\chi$ belongs to the kinetic space, and $D_+Q^*\chi=\Gamma_\HH\nabla\chi$ by closure. For every pair $(k,\nabla\chi)\in\cZ_+$ the difference $k-Q^*\chi$ is therefore an orbit invariant. This proves $\cZ_{\rm rem}=0$ in \eqref{eq:Schursep} and justifies writing $M_{z,\acc}=\mathcal L_z^0$ and $M_{\theta,\acc}=\mathcal L_\theta^0$ for these profiles.

\subsection{Constant axial drift: the Bennett equilibrium}
Write $\sigma_s$ for the species thermal speeds and $T_s=m_s\sigma_s^2$. In \cite[Theorem~1.1]{HolmerZhang2026},
\begin{equation}\label{eq:Bennett-profile}
 F_s=C_s\exp\big((-e_s+u_sp_{sz})/T_s\big),\qquad
 u_s=\operatorname{sgn}(q_s)\sigma_s,\qquad
 b=\frac{q_s u_s}{T_s}=\frac{|q_s|}{m_s\sigma_s}>0 .
\end{equation}
The two species have the same $b$. Choose the constants for pointwise neutrality and put $\phi^0=0$. Completing the velocity square gives mean $u_s$ and density $n_s=n_{s,c}\exp(b(A_z^0-A_c))$, where $A_c=A_z^0(0)$. Thus, with $j_c=\sum_s q_s u_s n_{s,c}$ and $c=b j_c/8$,
\begin{equation}\label{eq:Bennett-parameters}
 A_z^0=A_c-2b^{-1}\log(1+cr^2),\quad
 n_s=\frac{n_{s,c}}{(1+cr^2)^2},\quad
 K_z=b j_z^0=\frac{8c}{(1+cr^2)^2}.
\end{equation}
These formulas solve \eqref{eq:radial-Maxwell}, with $B_\theta^0$ decaying as $1/r$. Crews, Meier and Shumlak \cite{CrewsMeierShumlak2025} analyze the cyclotron and betatron orbits in this field. Sharma \cite{Sharma1983} used particle orbits to derive a kinetic eigenmode equation for the Bennett equilibrium and computed an electrostatic dispersion relation numerically. Meierovich \cite{Meierovich1986} studied $z$-independent Bennett disturbances near threshold in a relativistic isothermal two-fluid model. In the electromagnetic kinetic problem below a single vector of particle invariants suffices to determine the accessible index.

\begin{lemma}[One negative direction and one constraint]\label{lem:rankone}
Let $m$ be a closed, lower-bounded form of negative index one. If $m[u_R]\to Q_*<0$ and $m(a,u_R)\to\ip a f$ on a form core for some $f$ in the underlying Hilbert space, then $m[a]\ge |\ip a f|^2/Q_*$.
\end{lemma}
\begin{proof}
Complete the square on $\operatorname{span}\{a,u_R\}$. For large $R$, $m[u_R]<0$, so $m[a]-|m(a,u_R)|^2/m[u_R]\ge0$ by the index assumption. Pass to the limit and use form density.
\end{proof}

\begin{figure}[tb]
\centering\includegraphics[width=.77\linewidth]{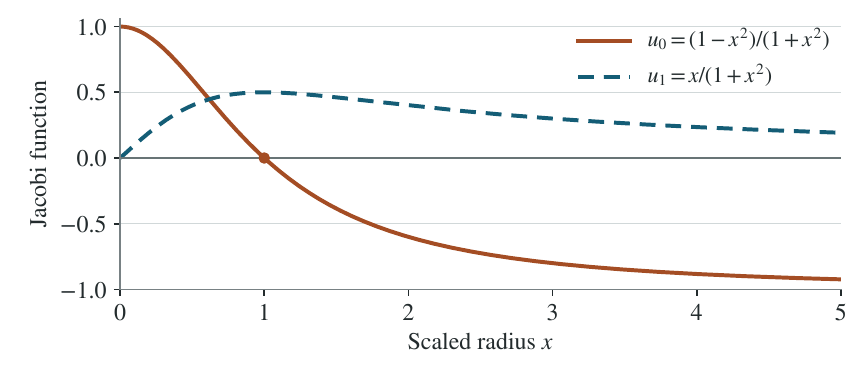}
\caption{Regular zero-energy solutions of the bare Bennett operator in scaled radius $x=\sqrt c\,r$. The radial solution $(1-x^2)/(1+x^2)$ changes sign once. The first angular solution $x/(1+x^2)$ is positive. Only the radial bare potential form has a negative direction, which the particle-number constraint removes.}\label{fig:bennett-jacobi}
\end{figure}

\begin{theorem}[Bennett pinch]\label{thm:Bennett-stable}
For \eqref{eq:Bennett-profile}--\eqref{eq:Bennett-parameters} on the whole transverse plane the bare index is one and $n^-(M_{z,\acc})=0$. The entire $z$-independent spectrum on the potential graph space is imaginary, including for unequal masses and temperatures satisfying \eqref{eq:Bennett-profile}.
\end{theorem}
\begin{proof}
The Maxwellians have $w_s=F_s/T_s$ and $\mu_{s,p_z}=u_sw_s$, so the coupling and shift bounds follow from the uniform zeroth and second velocity moments. Smooth magnetic fields conserve speed, giving a complete flow. In logarithmic radius Pfirsch \cite[\S2.5, Eqs.~(78)--(84)]{Pfirsch1962} obtained the $-2\sech^2$ potential and its marginal channels $j^2=0,1$. In the present radial coordinates the regular zero-energy solutions for $H_0,H_1$ are, up to constants,
\begin{equation}\label{eq:Bennett-Jacobi}
 y_0=\frac{1-cr^2}{1+cr^2},\qquad y_1=\frac{r}{1+cr^2}.
\end{equation}
The first has one zero. Sturm oscillation on disks with Dirichlet outer boundary, then exhaustion of the form core, gives radial negative index one. The second is positive for $r>0$, so the ground-state identity gives $H_1\ge0$. The potential difference $(j^2-1)/r^2$ gives $H_j\ge0$ for $|j|\ge1$. See Figure~\ref{fig:bennett-jacobi}.

For the accessible index set $\mathsf W=\sum_sT_s\int n_s\,\dd x_\perp$. The vector $((F_s)_s,0)$ belongs to $\cZ_{z,+}$ since each $F_s$ is constant along trajectories. As $F_s/w_s=T_s$, pairing $g$ with $(F_s)_s$ integrates each corrected density $g_s$ with weight $T_s$, and the squared norm of the vector is $\mathsf W$. The mean of $v_z$ is $u_s$, so the pairing of $\Gamma_z a$ with $(F_s)_s$ is $\ip{a}{j_z^0}$. Adding these two pairings cancels the density shift in $g$ and gives $\sum_s T_s\iint f_s\dd v\dd x_\perp$, a conserved linear combination of the species-number variations. Hence the correction in \eqref{eq:z-macc} is at least $|\ip a{j_z^0}|^2/\mathsf W$. On the other hand $j_z^0=b\sum_sT_s n_s$, so $K_z=b^2\sum_sT_s n_s$. The constant $-1/b$ solves $M_z(-1/b)=j_z^0$. Multiplying it by a smooth logarithmic cutoff equal to one on $r<R$ and zero beyond $R^2$ gives $u_R$ with $\|\nabla u_R\|_2^2=O(1/\log R)$. Consequently
\[
 m_z[u_R]\to-b^{-2}\int K_z\,\dd x_\perp=-\mathsf W,\qquad
 m_z(a,u_R)\to\ip a{j_z^0}.
\]
Lemma~\ref{lem:rankone} proves $m_z[a]\ge-|\ip a{j_z^0}|^2/\mathsf W$. The contribution of this one invariant cancels the lower bound. The remaining components of the projection contribute nonnegative squares. Apply Theorem~\ref{thm:z-application}.
\end{proof}

The same inequality has a short proof on the unit sphere. The conformal symmetry of the stationary Bennett profile also appears in the work of Kiessling and Lebowitz \cite{KiesslingLebowitz1994}. Stereographic projection in the variable $\sqrt{c}\,x_\perp$ gives $K_z\dd x_\perp=2\dd\omega$, where $\dd\omega$ is sphere area. Since $j_z^0=K_z/b$ and $\mathsf W=8\pi/b^2$, conformal invariance of the Dirichlet integral gives
\[
 m_z[a]+\frac{|\ip{a}{j_z^0}|^2}{\mathsf W}
   =\int_{S^2}|\nabla_{S^2}\widetilde a|^2\dd\omega
       -2\int_{S^2}|\widetilde a-\overline a|^2\dd\omega\ge0 .
\]
Here $\widetilde a$ is the pullback to the sphere and $\overline a=(4\pi)^{-1}\int_{S^2}\widetilde a\dd\omega$. The inequality is Poincar\'e's inequality with first positive spherical eigenvalue $2$. It holds first for compactly supported $a$ and then throughout the form domain by density.

The first angular solution in \eqref{eq:Bennett-Jacobi} also identifies a physical mode at the threshold. Infinitesimal translation in the $x_1$ direction gives the stationary variation
\[
 a=-\partial_1 A_z^0=\frac{4c}{b}\frac{r\cos\theta}{1+cr^2},
 \qquad f_s=-\partial_1F_s,\qquad E=0 .
\]
It has finite magnetic energy because $\nabla a\in L^2$. The density variation has finite kinetic norm because $f_s=bF_s a$ and $n_s=O(r^{-4})$. The potential itself has a $\cos\theta/r$ tail and fails to lie in $L^2(\R^2)$ for every choice of its additive constant. This physical neutral mode illustrates the role of the graph norm in Theorem~\ref{thm:Bennett-stable}. For a mode with $\lambda\ne0$, Faraday's law gives $a=-E_z/\lambda\in H^1(\R^2)$ since $\curl(a\eZ)$ is the transverse magnetic field.

\subsection{Variable axial drift}
Morse \cite[\S\S V--VI]{Morse1969} used a Gaussian dependence on quadratic canonical axial momentum in a different geometry. Channon and Coppins \cite{ChannonCoppins2001} constructed $z$-pinch equilibria with sheared axial flow in a Vlasov-fluid model. The two-species kinetic equilibrium in \cite[Theorem~1.2]{HolmerZhang2026} has
\begin{equation}\label{eq:z-variable-profile}
 F_s=C_s\exp\!\left(-e_s/T_s-\frac{p_{sz}^2}{2m_s^2\sigma_s^2}\right),
 \qquad \phi^0=0,\qquad b=|q_s|/(m_s\sigma_s).
\end{equation}
Now write $U=b A_z^0/\sqrt2$ and $n_{s,*}=C_s(2\pi\sigma_s^2)^{3/2}/\sqrt2$. Completing the square gives $n_s=n_{s,*}e^{-U^2/2}$ and $\langle v_z\rangle_s=-q_s A_z^0/(2m_s)$. The constants in \cite[Theorem~1.2]{HolmerZhang2026} give $\sum_s q_s n_{s,*}=0$. Hence, with $c_0=\frac12\sum_s q_s^2 n_{s,*}/m_s>0$ in our Maxwell units, Amp\`ere's equation and its current derivative become
\begin{equation}\label{eq:z-variable-current}
 -U''-r^{-1}U'=J(U),\qquad J(U)=-c_0 Ue^{-U^2/2},\qquad
 K_z=J'(U)=c_0(U^2-1)e^{-U^2/2}.
\end{equation}
Rescaling the radius multiplies $c_0$ by the square of the length scale. Remark~\ref{rem:normalizations} gives its value for the charge parameter of \cite{HolmerZhang2026}. Crews's isothermal flux-function model \cite{Crews2026} distinguishes a mean flow affine in flux from genuinely nonpolynomial flux dependences and exhibits an exponential example. In the family considered here the two species drift in opposite directions even at zero electric field, and their mean axial velocities are linear functions of $A_z^0$.

\begin{theorem}[Variable-drift pinch]\label{thm:z-variable}
For the regular nonconstant solutions of \eqref{eq:z-variable-current} the $z$-independent graph-space spectrum is contained in $\R\cup i\R$, and the restriction to all nonzero angular harmonics has spectrum in $i\R$. Thus any unstable $z$-independent mode belongs to the radial harmonic and has a real growth rate.
\end{theorem}
\begin{proof}
The Gaussian moments and density-shift coefficients are bounded since polynomials in $U$ times $e^{-U^2/2}$ are bounded. Smooth magnetic fields preserve speed. The right-hand side of the radial ODE is smooth and bounded, giving global continuation from $U'(0)=0$. By sign symmetry take $U(0)>0$. Then $(rU')'=c_0rUe^{-U^2/2}>0$ as long as $U>0$, so $U'>0$ throughout $r>0$. Differentiating the ODE gives $H_1U'=0$ with $U'=O(r)$ at the axis. For compactly supported smooth first-harmonic amplitudes the ground-state identity gives
\[
 h_1[a]=\int_0^\infty (U')^2
      \big|\partial_r(a/U')\big|^2 r\dd r\ge0.
\]
The intrinsic axis condition removes the boundary term. Form closure and \eqref{eq:z-Hj} give nonnegativity in every $|j|\ge1$ sector. Rotations reduce the full evolution and the accessible projection. The bare and accessible energies are nonnegative on the whole orthogonal sum of these sectors, so Theorem~\ref{thm:separable} applies on the sum. The sign of $M_{z,\acc}=\mathcal L_z^0$ on $j=0$ remains open. A search for growing modes can therefore use the radial potential equation \eqref{eq:growth-pencil}.
\end{proof}

\begin{remark}[Normalization of the four pinches]\label{rem:normalizations}
To compare the currents with \cite{HolmerZhang2026}, keep its charge parameter $a_{\mathrm{ch}}$, radial scale $r_0$, axial period $\zeta$, and species data fixed. The density prefactor for species $s$ has the same form in all four distributions,
\[
 C_s=\frac{a_{\mathrm{ch}}}{r_0^2\zeta\sigma_s^3|q_s|}.
\]
Write $D_0=(2\pi)^{3/2}(a_{\mathrm{ch}}/\zeta)\sum_s |q_s|/m_s$ in our Maxwell units. For the Bennett profile, $d_{\rm B}$ denotes the coefficient in $-\Delta_x(bA_z^0)=d_{\rm B}\exp(bA_z^0)$, where $x=x_\perp/r_0$. The four coefficients in the dimensionless field equations are
\[
 \begin{array}{ll}
 \text{Bennett axial drift} & d_{\rm B}=e^{1/2}D_0,\\
 \text{quadratic axial momentum, \eqref{eq:z-variable-current}} & r_0^2c_0=D_0/(2\sqrt{2}),\\
 \text{angular Gaussian, \eqref{eq:theta-G-ODE}} & d_\theta=D_0,\\
 \text{helical Gaussian, \eqref{eq:Gaussian-Volterra}} & r_0^2d_0=D_0 .
 \end{array}
\]
In SI units multiply $D_0$ by $\mu_0$. The constant axial drift gives the factor $\exp(u_s^2/(2\sigma_s^2))=e^{1/2}$ on completing the square. For the quadratic axial momentum the axial velocity integral contributes $1/\sqrt{2}$ and the mean is $-q_sA_z^0/(2m_s)$. The angular and helical density formulas retain their radius-dependent variance factors explicitly.

This distinguishes the coefficients denoted by $d$ in the introduction of \cite{HolmerZhang2026}. Equation (1.2) there gives $d_{\rm B}$, whereas the multiplier printed as $d/2$ in (1.10) is $r_0^2c_0$ for the distribution (1.12). In (1.5) the shifted axial square coming from (1.4) is $(v_z/\sigma^{\pm}\mp1)^2$, and (1.6) requires the additional factor $e^{1/2}$ with the stated charge parameter. Absorbing constants into a freely chosen density amplitude or radial scale leaves the profile shapes unchanged.
\end{remark}

\section{The \texorpdfstring{$\theta$}{theta}-pinch}\label{sec:theta}

Here an azimuthal current produces an axial magnetic field. Perturbations may depend on both $r$ and $z$ and preserve rotations. The family studied in \cite[Theorem~1.3]{HolmerZhang2026} has a Gaussian dependence on angular momentum. A second family has rigid rotation and an explicit density. We treat both on $\R^2\times\T$, with the axis understood in Cartesian coordinates.

\subsection{The axisymmetric operator}
Let $A^0=A_\theta^0(r)\eTh$, $B_z^0=(rA_\theta^0)'/r$, $E^0=-\phi^{0\prime}\eR$ and $F_s=\mu_s(e_s,\ell_s)$ with $w_s=-\mu_{s,e}>0$. Amp\`ere's equation is $-(B_z^0)'=j_\theta^0$. The symmetry reflection changes $(v_r,v_z)$ to $(-v_r,-v_z)$ and fixes $v_\theta$. To see reversal in the reduced particle equations, observe that the $v_r$ acceleration is $v_\theta^2/r+(q_s/m_s)(E_r^0+v_\theta B_z^0)$, the $v_\theta$ acceleration is $-v_r(v_\theta/r+q_sB_z^0/m_s)$, and $v_z$ is constant. Rotation has zero twist and $B_\theta^0=0$.

Writing $a=A_\theta$, the momentum correction and current coefficient become
\[
 g_s=f_{s,\ev}-q_s r\mu_{s,\ell}a,\qquad
 (\Gamma_\theta a)_s=q_sw_sv_\theta a,\qquad
 K_\theta=\sum_s q_s^2r\int v_\theta\mu_{s,\ell}\dd v=\partial_{A_\theta}j_\theta^0 .
\]
The curl of $a\eTh$ has components $-a_z\eR+(a_r+a/r)\eZ$, so
\begin{align}
 m_\theta[a]&=\int\left(|a_r|^2+r^{-2}|a|^2+|a_z|^2-K_\theta|a|^2\right)r\dd r\dd z, \label{eq:theta-form}\\
 \ip{M_{\theta,\acc} a}{a}_{\rm form}
 &=m_\theta[a]+\|P_{\cZ_{\theta,+}}(\Gamma_\theta a,0)\|^2 .\label{eq:theta-macc}
\end{align}
Here $\cZ_{\theta,+}$ is defined as in \eqref{eq:z-macc} with $E_\perp=(E_r,E_z)$ and velocities $(v_r,v_z)$. Smooth core amplitudes are $O(r)$ at the axis. Integration by parts on this core accounts for the $r^{-2}$ term. The closed forms describe the natural graph space of $a\eTh$.

\begin{theorem}[Axisymmetric spectrum]\label{thm:theta-application}
For a smooth stationary $\theta$-pinch as above, assume the flow, moment and density-shift hypotheses of Section~\ref{sec:framework}. The entire axisymmetric spectrum on the potential graph space lies in $\R\cup i\R$, and it lies in $i\R$ precisely when $M_{\theta,\acc}\ge0$. For finite index the positive-real algebraic count equals $n^-(M_{\theta,\acc})$. The assertions hold on the full, constrained and accessible energy spaces, also on a circular conducting cylinder with specular particles and $a(R,z)=0$.
\end{theorem}
\begin{proof}
The reflection, curl and current above give the separable system of Theorem~\ref{thm:separable}, with accessible energy \eqref{eq:theta-macc}. The specular wall map commutes with the reflection. On a finite axially periodic cylinder the curl-free $(E_r,E_z)$ fields in the electric domain are gradients. Indeed, the axial period of $E_z$ is independent of $r$ and vanishes at the wall. Thus the scalar electrostatic solve in Section~\ref{sec:framework} exhausts the electric correction on that cell.
\end{proof}
For related sharp criteria with conducting walls and specular particles in the relativistic system, see Nguyen and Strauss \cite{NguyenStrauss2013,NguyenStrauss2014} and Zhang \cite{Zhang2019} for general axisymmetric domains.

\subsection{Gaussian angular momentum}
For \cite[Theorem~1.3]{HolmerZhang2026}, let $r_0$ be the radial scale, $x=r/r_0$ and $U(x)=b A_\theta^0(r)$ with $b=|q_s|/(m_s\sigma_s)$ common to the two species. Then
\begin{equation}\label{eq:theta-gaussian}
 F_s=C_s\exp\left(-e_s/T_s-\frac{\ell_s^2}{2m_s^2\sigma_s^2r_0^2}\right),\qquad \phi^0=0 .
\end{equation}
Set $n_{s,0}=C_s(2\pi\sigma_s^2)^{3/2}$ with $\sum_s q_s n_{s,0}=0$. The angular Gaussian integral gives
\[
 t=\frac{x^2 U^2}{1+x^2},\qquad
 n_s=n_{s,0}(1+x^2)^{-1/2}e^{-t/2},\qquad
 \langle v_\theta\rangle_s=-\frac{q_s A_\theta^0}{m_s}\frac{x^2}{1+x^2}.
\]
Consequently the dimensionless current and its derivative are
\begin{align}
 -U''-x^{-1}U'+x^{-2}U
   &=J_\theta(x,U)=-d_\theta\frac{x^2 U}{(1+x^2)^{3/2}}e^{-t/2},\label{eq:theta-G-ODE}\\
 r_0^2 K_\theta(r)&=\partial_UJ_\theta(x,U)
    =d_\theta\frac{x^2}{(1+x^2)^{3/2}}(t-1)e^{-t/2}.\label{eq:theta-G-K}
\end{align}
Here $d_\theta=r_0^2\sum_s q_s^2 n_{s,0}/m_s$ follows from $br_0^2j_\theta^0=J_\theta$. Remark~\ref{rem:normalizations} gives the coefficient in \cite[\S6.1, Eq.~(6.5)]{HolmerZhang2026} with the charge parameter held fixed.

\begin{theorem}[A stable range for the Gaussian $\theta$-pinch]\label{thm:theta-gaussian}
Let $s=|U'(0)|>0$ on the regular branch $U(0)=0$ of \eqref{eq:theta-G-ODE}. If
\[
 s\ge C_*d_\theta,\qquad
 C_*=\sqrt{2+\sqrt3}\,(1+\sqrt3)\exp[-(2+\sqrt3)/2]<0.817,
\]
the axisymmetric spectrum on the potential graph space is imaginary. In particular this holds at $s=d_\theta=1$. For every $s,d_\theta>0$ the axes conclusion of Theorem~\ref{thm:theta-application} holds.
\end{theorem}
\begin{proof}
Choose the sign of $U$ so that $U'(0)=s$. Write $U=x h$. Equation~\eqref{eq:theta-G-ODE} gives
\[
 (x^3h')'=d_\theta x^5(1+x^2)^{-3/2} e^{-t/2}h,\qquad h(0)=s .
\]
It follows by integration that $h\ge s$ and $U\ge sx$ on the whole branch. The bounded coefficient in the corresponding linear equation for $h$ on compact intervals gives global continuation. In particular $t\ge s^2x^4/(1+x^2)$, so the density has a Gaussian tail. Completing the velocity square also bounds all couplings and the density-shift moment $r^2\int |\mu_{s,\ell}|^2/w_s\,\dd v$: the shifted Gaussian has angular variance $\sigma_s^2/(1+x^2)$, and powers of $U$ in these moments are controlled by $e^{-t/2}$ for $x\ge1$. The complete flow follows from smoothness and conservation of particle speed.

Multiplying \eqref{eq:theta-G-K} by $x^2$, use
\[
 \frac{x^4}{(1+x^2)^{3/2}}
 \le \frac{x^2}{\sqrt{1+x^2}}\le \frac{\sqrt{t}}s .
\]
The maximum of $\sqrt t(t-1)_+e^{-t/2}$ occurs at $t=2+\sqrt3$. Hence $r^2(K_\theta)_+\le C_*d_\theta/s$. The angular term $r^{-2}$ in \eqref{eq:theta-form} dominates the positive part of $K_\theta$ throughout the stated range. Apply \eqref{eq:theta-macc} and Theorem~\ref{thm:theta-application}.
\end{proof}

\subsection{Rigid rotation and an explicit positive solution}
Komarov and Fadeev \cite[Eqs.~(4)--(6)]{KomarovFadeev1962} obtained neutral multicomponent rigid rotors with density proportional to $\sech^2(ar^2+b)$ and an axial field affine in $\tanh(ar^2+b)$. Dickman, Morse and Nielson \cite[Eqs.~(8)--(9)]{DickmanMorseNielson1968} gave the profile and its cylindrical marginal equation. Freidberg and Morse \cite{FreidbergMorse1969} related collisionless tearing stability to the position of the density maximum. For mirror-type modes Schindler \cite[p.~739]{Schindler1966} deduced stability of the cylindrical examples treated by Pfirsch \cite{Pfirsch1962}. Batchelor and Davidson \cite{BatchelorDavidson1975} studied kinetic rigid rotors using the energy in a rotating frame. Crews \cite{Crews2026} derived the Liouville equation for the flux in the $\theta$-pinch limit of the cylindrical construction of Cordonnier et al.\ \cite{Cordonnier2022}. Here we evaluate the full axisymmetric graph-space spectrum for the following center-peaked, neutral counter-rotating member.

\begin{figure}[tb]
\centering\includegraphics[width=.81\linewidth]{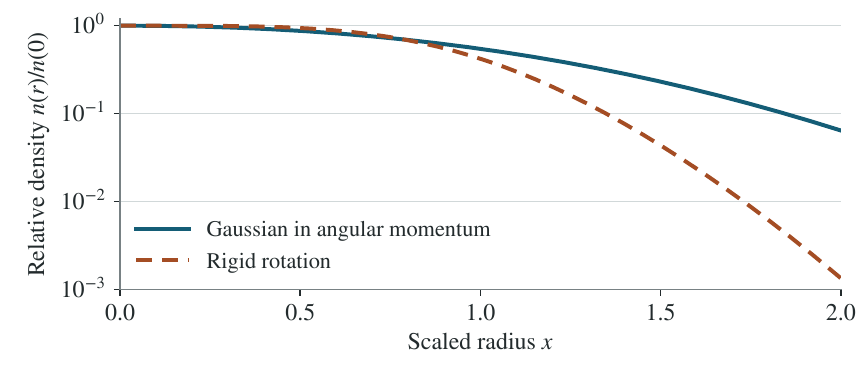}
\caption{Density relative to its on-axis value. The Gaussian curve comes from integrating the equilibrium equation, with $x=r/r_0$ and $s=d_\theta=1$ in \eqref{eq:theta-G-ODE}. For the rigid rotor $x=\sqrt{k}\,r$ and $n/n_c=\sech^2(x^2)$. At scaled radius $x=2$, the curves reach $6.42\%$ and $0.134\%$ of their respective central densities.}\label{fig:theta-comparison}
\end{figure}

Take charges $q_\pm=\pm q$, equal mass $m$, and constants $\beta,\Omega,n_c>0$. Put $\alpha=\beta q\Omega$, $k=q\Omega\sqrt{\beta n_c}/2$, and
\begin{equation}\label{eq:theta-constants}
\begin{aligned}
 F_\pm&=n_c(\beta m/2\pi)^{3/2}\exp(-\beta e_\pm\pm\beta\Omega\ell_\pm),\\
 A_\theta^0&=-\frac{m\Omega}{2q}r-\frac{2\log\cosh(kr^2)}{\alpha r},\qquad \phi^0=0 .
\end{aligned}
\end{equation}
The species have angular means $\pm\Omega r$ and the same density $n(r)=n_c\sech^2(kr^2)$. Figure~\ref{fig:theta-comparison} compares this density with the angular Gaussian. The current, magnetic field and current derivative are
\[
 j_\theta^0=2q\Omega r n(r),\qquad
 B_z^0=-\frac{m\Omega}{q}-\frac{4k}{\alpha}\tanh(kr^2),\qquad
 K_\theta=8k^2r^2\sech^2(kr^2).
\]
Differentiation verifies $-(B_z^0)'=j_\theta^0$. Smoothness, conservation of speed and the fast decay of the shifted Gaussian moments give the hypotheses of Theorem~\ref{thm:theta-application}.

\begin{theorem}[Rigidly rotating $\theta$-pinch]\label{thm:theta-stable}
For \eqref{eq:theta-constants}, the entire axisymmetric spectrum on the potential graph space lies in $i\R$.
\end{theorem}
\begin{proof}
The regular marginal solution found by Pfirsch \cite[\S2.4, Eqs.~(63)--(68)]{Pfirsch1962} is $\varphi(r)=\tanh(kr^2)/r$ in these variables. This positive function solves
$(-\partial_r^2-r^{-1}\partial_r+r^{-2}-K_\theta)\varphi=0$ by differentiation. It behaves as $kr$ at the axis. Integration by parts on the core gives
\[
 m_\theta[a]=\int\left(\varphi^2
       \left|\partial_r(a/\varphi)\right|^2+|\partial_z a|^2\right) r\dd r\dd z\ge0 .
\]
Taking the form closure proves the assertion by Theorem~\ref{thm:theta-application}. On the whole radial half-line $\varphi\sim r^{-1}$ is a threshold solution. Its positivity suffices for the form identity.
\end{proof}

\section{The screw-pinch}\label{sec:screw}

In the screw-pinch of \cite[Theorem~1.4]{HolmerZhang2026}, the current at radius $r$ follows the helix tangent to $X_\kappa=\partial_\theta+\kappa\partial_z$, where $\kappa\ne0$ is its axial advance per radian. The magnetic field lies in the perpendicular direction on the same cylinder. In the orthonormal frame
\begin{equation}\label{eq:helical-frame}
 \delta=(r^2+\kappa^2)^{1/2},\qquad
 \eeta=(r\eTh+\kappa\eZ)/\delta,\qquad
 \esig=(-\kappa\eTh+r\eZ)/\delta ,
\end{equation}
write $\alpha=A_\eta^0$. The particle momentum along the helix and the equilibrium distribution are
\begin{equation}\label{eq:Gaussian-profile}
 P_s=\ell_s+\kappa p_{sz}=\delta(m_sv_\eta+q_s\alpha),\qquad
 F_s=c_s\exp(-\beta_se_s-\gamma_sP_s^2/2),\quad \phi^0=0 .
\end{equation}
To compare with \cite[Eq.~(1.25)]{HolmerZhang2026}, write $r_0$ for the length denoted there by $\beta$, $\zeta$ for the axial period, $\sigma_s$ for the thermal speeds, and $a_{\rm ch}$ for the charge parameter. In the notation of \eqref{eq:Gaussian-profile},
\[
 \beta_s=(m_s\sigma_s^2)^{-1},\qquad
 \gamma_s=(m_s^2\sigma_s^2r_0^2)^{-1},\qquad
 c_s=\frac{a_{\rm ch}}{r_0^2\zeta\sigma_s^3|q_s|}.
\]
The two species have a common value of $|q_s|/(m_s\sigma_s)$. Morse used separate axial and angular momentum multipliers to construct equilibria with drift and rotation \cite[\S\S V--VI]{Morse1969}. Virtamo combined translation and rotation in locally Maxwellian equilibria and computed screw-pinch profiles \cite{Virtamo1980}. Cordonnier, Leoncini, Dif-Pradalier and Garbet constructed a self-consistent screw family with separate linear multipliers and rigid relative drifts \cite{Cordonnier2022}. Cordonnier, Lebouazda, Leoncini and Dif-Pradalier subsequently studied particle motion in these fields \cite{Cordonnier2026}. Crews treats sheared equilibria with flow depending on the flux \cite{Crews2026}. For stability, Mondt and Davidson derived a low-frequency screw-pinch eigenvalue equation with Vlasov ions and drift-kinetic electrons \cite{MondtDavidson1980}. In \eqref{eq:Gaussian-profile} the square of $\ell_s+\kappa p_{sz}$ couples the two momenta. This single conserved combination determines which perturbations admit the energy reduction for the two-species Vlasov--Maxwell equations.

Throughout this section the parallel potential belongs to $D(C_e)$ with the graph norm \eqref{eq:Lambda}. The full, constrained and accessible spectra refer to the realizations on this space defined in Section~\ref{sub:closed}. The requirement that the potential be square integrable matters on the whole cylinder.

\subsection{Radial perturbations}\label{sub:radial-screw}
First allow any smooth radial equilibrium $F_s=\mu_s(e_s,\ell_s,p_{sz})$ with radial electric field and tangential magnetic field. A radial perturbation retains both momenta. Its two tangential potential components are $\mathbf a=(a_\theta,a_z)^T$. The shift and current response of Section~\ref{sec:framework} become
\begin{align}
 g_s&=f_{s,\ev}-q_s(r\mu_{s,\ell}a_\theta+\mu_{s,p_z}a_z),&
 (\Gamma_\VV^r\mathbf a)_s&=q_sw_s(v_\theta a_\theta+v_za_z),\nonumber\\
 K_r&=\sum_s q_s^2\int \binom{v_\theta}{v_z}
                    (r\mu_{s,\ell},\mu_{s,p_z})\dd v,&
 C_r\mathbf a&=(-a_z',a_\theta'+a_\theta/r)^T .
 \label{eq:screw-radial-K}
\end{align}
Here $w_s=-\mu_{s,e}$ and $K_r$ is symmetric by velocity integration by parts. The radial measure is $r\dd r$, with constant orbit-volume factors suppressed. Smooth axis data have $a_\theta=r a_1(r^2)$, $a_z=a_2(r^2)$ near zero. Close the curl on this intrinsic core, imposing $\mathbf a(R)=0$ at a conducting wall. For the whole cylinder take $R=\infty$.

At fixed momenta $(\ell,p_z)$, the radial characteristic is the mechanical system
\begin{equation}\label{eq:effective-U}
 D_s^r=v_r\partial_r-m_s^{-1}U_s'\partial_{v_r},\qquad
 U_s=\frac{(\ell-q_srA_\theta^0)^2}{2m_sr^2}
             +\frac{(p_z-q_s A_z^0)^2}{2m_s}+q_s\phi^0 .
\end{equation}
Reflection of $v_r$ reverses this flow. Write $D_+^r,D_-^r$ for its restrictions to even and odd inputs. With $(\Gamma_rE_r)_s=q_sw_sv_rE_r$, put
\[
 \cZ_r=\{(g,E_r):g\in D(D_+^r),\ D_+^rg=\Gamma_rE_r\},\qquad
 \mathfrak l_r[\mathbf a]=\norm{C_r\mathbf a}^2-\ip{K_r\mathbf a}{\mathbf a}
                  +\norm{P_{\cZ_r}(\Gamma_\VV^r\mathbf a,0)}^2 .
\]
Let $\mathcal L^0_{\rm rad}$ denote the self-adjoint operator of this closed form. It acts on two spatial components. On a finite solid cell $E_r=\chi'$ with $\chi(R)=0$, so the electrostatic equation can be solved as in Section~\ref{sec:framework}. The invariant particle projection is the time average along each connected radial orbit of \eqref{eq:effective-U}. On a periodic trajectory this is the normalized integral with weight $\dd r/|v_r|$ on each leg. Only $\ell=0$ trajectories can reach the axis, and they form a set of phase volume zero. The Cartesian flow fixes the transport domain there.

\begin{theorem}[Radial screw-pinch spectrum]\label{thm:radial-screw}
Under the flow, moment and shift hypotheses of Section~\ref{sec:framework}, the entire radial spectrum lies in $\R\cup i\R$. It is imaginary if and only if $\mathcal L^0_{\rm rad}\ge0$. When $N_r=\nm(\mathcal L^0_{\rm rad})<\infty$, the positive real eigenvalues have total algebraic multiplicity $N_r$. These assertions hold on the full radial space and on the constrained and accessible spaces.
\end{theorem}
\begin{proof}
The curl, current matrix and shift are given by \eqref{eq:screw-radial-K}. At a specular wall reflection reverses the direction of traversal in \eqref{eq:effective-U}. Thus there is strong radial reversal on the closed transport domain. With $\Lambda_r=(I+C_r^*C_r)^{1/2}$, the off-diagonal coupling of Section~\ref{sec:framework} is
\[
 (h,\mathbf E_\VV)\longmapsto
 (-D_-^rh+\Gamma_\VV^r\mathbf E_\VV,-\Gamma_r^*h,-\Lambda_r\mathbf E_\VV)
\]
on $D(D_-^r)\oplus D(\Lambda_r)$. Theorems~\ref{thm:accessibility} and \ref{thm:separable} give the form $\mathfrak l_r$ and the stated count.
\end{proof}

\subsection{Perturbations constant along a helix}\label{sub:helical-application}
A helical perturbation may depend on both $r$ and $\sigma=z-\kappa\theta$. For this class the retained canonical momentum is $P_s$. A profile $\mu_s(e,\ell,p_z)$ factors through $(e,P_s)$ precisely when $\mu_{s,p_z}=\kappa\mu_{s,\ell}$ on connected momentum fibers. For \eqref{eq:Gaussian-profile} write the factor as $\nu_s(e,P_s)$. In \cite[\S7]{HolmerZhang2026} the axial period is $\zeta$ and $2\pi\kappa/\zeta\in\mathbb Z$, so $\sigma$ has period $L=\zeta$. More generally, the angular and axial identifications give the two periods $2\pi\kappa$ and $\zeta$ for $\sigma$. If $2\pi\kappa/\zeta=p/q$ in lowest terms, use $L=\zeta/q$ and $n_\kappa=2\pi\kappa/L=p$. If $2\pi\kappa/\zeta$ is irrational, the helical subgroup is dense in the rotation--translation torus. Strong continuity of the lifted torus action then makes every invariant $L^2$ state invariant under both rotations and translations, so only the radial sector remains. In the Fourier discussion below we use $L$ for the rational case. In Fourier harmonic $j$, a smooth invariant scalar has angular order $|n_\kappa j|$ at the axis. For vector and kinetic fields we use the intrinsic Cartesian core of Section~\ref{sub:closed}.

The helical reflection changes $(v_r,v_\sigma)$ to $(-v_r,-v_\sigma)$ and leaves $v_\eta$ fixed. The frame itself rotates along a trajectory. Differentiating \eqref{eq:helical-frame} in the Lorentz transport gives the part that preserves parity,
\begin{equation}\label{eq:helical-Dpre}
 D_{s,\pre}=\left(\tau v_\eta+\frac{q_s}{m_s}B_\eta^0\right)
              (v_r\partial_{v_\sigma}-v_\sigma\partial_{v_r}),
 \qquad \tau=\frac{2\kappa}{\delta^2}.
\end{equation}
For the Gaussian equilibrium \eqref{eq:Gaussian-profile} considered here, $B_\eta^0=0$, so the term proportional to $\tau v_\eta$ persists. The same $\tau$ appears in the curl. If $a=A_\eta$, the perpendicular magnetic perturbation from $a$ is
\[
 Ca=-\frac\delta r a_\sigma\eR+
             \left(a_r+\frac r{\delta^2}a\right)\esig .
\]
Write $b=B_\eta$ for the total parallel magnetic perturbation. Then $b-\tau a$ is the parallel curl of the perpendicular potential. We use the electric and magnetic realizations $C_e,C_m$ of $C$ from Section~\ref{sec:framework}. At a conducting wall these are respectively the minimal and maximal realizations. On a finite solid cell the flux condition is $\ip{b-\tau a}{\delta^{-1}}=0$ since $\Ker C_m=\spanop\{\delta^{-1}\}$.

Let $(\Gamma_i u)_s=q_sw_sv_i u$ for $i=r,\eta,\sigma$, and $\Gamma_\HH E_\HH=\Gamma_r E_r+\Gamma_\sigma E_\sigma$. The shifted density is $g_s=f_{s,\ev}-q_s\delta\nu_{s,P}a$. Put $K_\eta=\sum_s q_s^2\delta\int v_\eta\nu_{s,P}\dd v$. With measure $r\dd r\dd\sigma$, the bare and accessible forms are
\begin{align}
 m_\kappa[a]&=\norm{C_ea}^2-\ip{K_\eta a}{a},\nonumber\\
 m_{\kappa,\acc}[a]&=m_\kappa[a]+
       \norm{P_{\cZ_\kappa}(\Gamma_\eta a,0,-\tau a)}^2,\label{eq:helical-macc}\\
 \cZ_\kappa&=\{(k,E_\HH,b):Dk=\Gamma_\HH E_\HH,\
       C_m b=\Gamma_\HH^*P_-k,\ C_m^*E_\HH=0\}.\nonumber
\end{align}
Here $D$ is the full characteristic generator, $P_-$ extracts the odd part under the helical reflection, and every equation defining $\cZ_\kappa$ uses the indicated operator domain. The joint projection retains the magnetic flux and particle invariants in a single energy.

\begin{theorem}[Helical screw-pinch spectrum]\label{thm:helical-application}
For a stationary profile $\nu_s(e_s,P_s)$ satisfying the hypotheses of Section~\ref{sec:framework}, nonnegativity of $m_{\kappa,\acc}$ places the entire helical spectrum on $i\R$. At finite index $N_\kappa=\nm(m_{\kappa,\acc})$, the nonimaginary spectrum consists of finitely many eigenvalues of finite algebraic multiplicity and $k_r+2k_c\le N_\kappa$. Here $k_r$ counts positive real eigenvalues and $k_c$ counts those in the open first quadrant, with algebraic multiplicity. Thus $N_\kappa\le1$ places the spectrum on the two axes.
\end{theorem}
\begin{proof}
The helical shift, curl and parity calculation identify the entries of the general Hamiltonian system in Section~\ref{sec:framework}. The term $a\mapsto-\tau a$ in the projection is bounded because $|\tau|\le2/|\kappa|$. Theorem~\ref{thm:accessibility} gives \eqref{eq:helical-macc}. Theorem~\ref{thm:general-spectral} gives the bound and the Hamiltonian reflection and conjugation symmetries on the full, constrained and accessible spaces.
\end{proof}

\paragraph{Computing the projection on a finite cylinder.}
On $0<r<R$ with a conducting wall, the scalar solve reveals how the twist changes the energy. Write $\zeta_m=\delta^{-1}$ for the magnetic null vector. If $C_m^*E_\HH=0$, the one-form perpendicular to the helix, $E_r\dd r+(r/\delta)E_\sigma\dd\sigma$, is closed and has vanishing boundary tangential trace. Lifting to the disk, taking a primitive, and fixing the boundary constant gives $E_\HH=\nabla\chi$ with $\chi\in H_0^1$. The primitive descends under the cyclic identification by averaging. For an invariant scalar $b$, \eqref{eq:curl-one} gives $C_m b=\delta^{-1}\nabla(\delta b)\times\eeta$. The gradient is perpendicular to $\eeta$. Since $\delta$ is smooth, bounded and bounded away from zero on this cell, the maximal magnetic graph norm is equivalent to the $H^1$ norm. Rellich compactness then gives $\norm b\le c_R\norm{C_m b}$ for $b\perp\zeta_m$. Indeed, a normalized sequence contradicting this estimate would converge along a subsequence to a unit vector in $\Ker C_m\cap\zeta_m^\perp$. Thus $\Ran C_m=(\Ker C_m^*)^\perp$ is closed, and $C_m^\dagger$ denotes its inverse with values perpendicular to $\zeta_m$. For $\ell\in\Ker D$, one has $\ip{\Gamma_\HH^*P_-\ell}{\nabla\chi}=\ip{\ell}{DQ^*\chi}=0$. Thus we may set
\[
 \Pi=P_{\Ker D},\qquad T\ell=C_m^\dagger\Gamma_\HH^*P_-\ell,\qquad
 H=I+T^*T,\qquad r_a=\Pi\Gamma_\eta a-T^*(\tau a).
\]
The equations for $\cZ_\kappa$ give its elements as
$(Q^*\chi+\ell,\nabla\chi,T\ell+c\zeta_m)$.
Minimizing the distance from $(\Gamma_\eta a,0,-\tau a)$ to these vectors first over $c$ and $\ell$ leaves one scalar variational equation. Define
\begin{align*}
 d(\chi,\psi)&=\ip{\nabla\chi}{\nabla\psi}+\ip{Q^*\chi}{Q^*\psi}
                     -\ip{H^{-1}\Pi Q^*\chi}{\Pi Q^*\psi},\\
 d(\chi_a,\psi)&=\ip{\Gamma_\eta a}{Q^*\psi}
                     -\ip{H^{-1}r_a}{\Pi Q^*\psi}\qquad(\psi\in H_0^1).
\end{align*}
Because $0<H^{-1}\le I$, $d(\chi,\chi)\ge\norm{\nabla\chi}^2$ and $\chi_a$ exists uniquely. Completing the squares gives
\begin{equation}\label{eq:helical-Schur}
 m_{\kappa,\acc}[a]=m_\kappa[a]+\ip{H^{-1}r_a}{r_a}
                 +d(\chi_a,\chi_a)+\frac{|\ip{\tau a}{\zeta_m}|^2}{\norm{\zeta_m}^2}.
\end{equation}
The last term is the cost of fixing the magnetic flux. The first two corrections use the full helical particle orbits. On the whole cylinder $\zeta_m\notin L^2(r\dd r)$, and \eqref{eq:helical-macc} supplies the projection without a closed-range assumption.

\subsection{The equilibrium equation and its coefficients}\label{sub:Gaussian}
For the two species of \eqref{eq:Gaussian-profile}, the parameters $a_*$ and $b_*$ below are independent of $s$. We assume $m_s,\beta_s,\gamma_s,c_s>0$ and $q_s\ne0$, and allow any family with these common parameters and the stated neutrality condition. In particular $a_*,b_*,d_0$ are positive. Write
\begin{equation}\label{eq:Gaussian-neutrality}
 a_*=\gamma_sm_s/\beta_s,\quad b_*=\gamma_sq_s^2,\quad
 n_{s,0}=c_s(2\pi/(\beta_sm_s))^{3/2},\quad
 \sum_s q_sn_{s,0}=0,\quad d_0=\sum_s q_s^2n_{s,0}/m_s .
\end{equation}
With the normalization of \cite[Eq.~(1.25)]{HolmerZhang2026}, $a_*=r_0^{-2}$ and $b_*=(|q_s|/(m_s\sigma_sr_0))^2$. Completing the velocity square in \eqref{eq:Gaussian-profile} gives
\begin{align}
 h&=1+a_*\delta^2,\qquad t=b_*\delta^2\alpha^2/h,\qquad
 V=d_0(h-1)h^{-3/2}e^{-t/2},\nonumber\\
 n_s&=n_{s,0}h^{-1/2}e^{-t/2},\qquad
 j_\eta^0=-V\alpha,\qquad K_\eta=V(t-1).\label{eq:Gaussian-current}
\end{align}
Indeed $v_\eta$ has mean $-\gamma_s q_s\delta^2\alpha/(\beta_sh)$ and variance $(\beta_sm_sh)^{-1}$, while $v_r,v_\sigma$ have zero mean. Then $K_\eta=\partial_\alpha j_\eta^0$ at fixed $r$. For radial perturbations, writing $\mathbf t=(r,\kappa)^T/\delta$ gives $K_r=K_\eta\mathbf t\mathbf t^T$. Thus both perturbation classes use the same scalar current response.

For $B_\eta^0=0$, the stationary Maxwell equations take a particularly simple form with $y=\delta\alpha$ and $p=r/\delta^2$,
\begin{equation}\label{eq:Gaussian-Volterra}
 (p y')'=pV(r,y/\delta)y,\qquad y(0)=|\kappa|\alpha_0,\quad y'(0)=0.
\end{equation}
The remaining potential obeys $(r A_\sigma^0/\delta)'=2\kappa r\alpha/\delta^3$ with zero integration constant. Then $B_\sigma^0=y'/\delta$, and $B_r^0=E^0=0$. With these species parameters and radial scale $r_0$, \eqref{eq:Gaussian-Volterra} has the form of \cite[Eq.~(1.24)]{HolmerZhang2026}. In our Maxwell units its coefficient is $r_0^2d_0=(2\pi)^{3/2}(a_{\rm ch}/\zeta)\sum_s|q_s|/m_s$, as obtained from the current in \cite[\S7.1, Eqs.~(7.12)--(7.13)]{HolmerZhang2026}. Remark~\ref{rem:normalizations} compares this value with the constant-drift coefficient for the same input normalization.

\begin{proposition}[The whole-cylinder Gaussian branch]\label{prop:Gaussian-global}
Under the positive-parameter and neutrality hypotheses of \eqref{eq:Gaussian-neutrality}, the equilibrium is smooth on the whole cylinder for each $\alpha_0>0$. For some $I_\infty>0$,
$\alpha(r)=I_\infty r/2+O(r^{-1}\log r)$ and $B_\sigma^0(r)\to I_\infty$. Both perturbation classes satisfy the flow, moment and shift hypotheses.
\end{proposition}
\begin{proof}
Integration of \eqref{eq:Gaussian-Volterra} gives $(py')(r)=\int_0^rpVy\dd s$. The resulting Volterra equation for $y$ is locally contractive at the regular axis and gives a smooth even solution. The positive initial value makes $y$ increasing. At most linear growth of the right-hand side in $y$ on compact radial intervals gives continuation. For any $r_1>0$, $(py')(r)\ge (py')(r_1)>0$ for $r\ge r_1$. Integrating $p^{-1}=r+\kappa^2/r$ now gives $y(r)\gtrsim r^2$ and hence $t=b_*y^2/h\gtrsim r^2$ at infinity. Since $y=\sqrt{ht/b_*}$, the factor $y e^{-t/2}$ is bounded by a constant times $\sqrt h\,e^{-t/4}$. The integrand $pVy$ therefore has a Gaussian tail. Thus $py'=I_\infty+O(e^{-cr^2})$ for some $I_\infty,c>0$. Integrating $p^{-1}=r+\kappa^2/r$ gives the asserted asymptotics. The moment and squared shift coefficients are bounded by polynomials in $r,\alpha$ times $h^{-1/2}e^{-t/2}$ by Gaussian integration. They decay at infinity and are bounded near zero. The Cartesian field is smooth, and $E^0=0$ preserves speed, giving a complete volume- and weight-preserving characteristic flow.
\end{proof}

\subsection{Which helical harmonics can grow?}\label{sub:high}
A radial equilibrium commutes with translations in $\sigma$. The Fourier projections reduce the transport group, Maxwell operators and projection in \eqref{eq:helical-macc}, with frequency $\nu_j=2\pi j/L$. In the zero harmonic, only $r$ varies. We may therefore change the particle parity back to reflection of $v_r$ alone. This gives a separable system within the original helical space.

\begin{theorem}[Zero helical harmonic]\label{thm:zero}
For a radial equilibrium as in Theorem~\ref{thm:helical-application}, the spectrum in harmonic zero lies in $\R\cup i\R$. If $N_0=\nm(m_{\kappa,\acc,0})<\infty$, there are exactly $N_0$ positive real eigenvalues counted with algebraic multiplicity.
\end{theorem}
\begin{proof}
Split the corrected density into its radial even and odd parts $g^r,h^r$, and write $T_mb=b'+rb/\delta^2$ on the magnetic domain, $\Lambda_0=(I+C_{e,0}^*C_{e,0})^{1/2}$. Regroup the variables into
$U=(g^r,E_r,\Lambda_0 a,b)$ and $V=(h^r,E_\eta,E_\sigma)$. The off-diagonal coupling from $V$ to $U$ is
\[
 \cB_0(h,e_\eta,e_\sigma)=
 \big(-D_-^rh+\Gamma_\eta e_\eta+\Gamma_\sigma e_\sigma,
       -\Gamma_r^*h,-\Lambda_0e_\eta,-\tau e_\eta-T_m^*e_\sigma\big)
\]
on $D(D_-^r)\oplus D(\Lambda_0)\oplus D(T_m^*)$. This follows from the radial transport reversal, $(E_\sigma)_t=-\Gamma_\sigma^*g^r+T_mb$, and $b_t=-\tau E_\eta-T_m^*E_\sigma$. The derivative entries act on separate components, so their direct sum plus the bounded entries is closed. In these coordinates the $V$ energy is the identity, and only the potential component of the $U$ energy differs from it. The self-adjoint operator for the squared spectral parameter in Theorem~\ref{thm:separable} is $\cB_0^*\cB_0-\Pi_\eta^*(I+K_\eta)\Pi_\eta$, where $\Pi_\eta V=E_\eta$. The same skew operator in the regrouped coordinates has off-diagonal blocks $\cB_0,-\cB_0^*$. Thus the closure of its range, which defines accessibility, is unchanged under the permutation and radial parity decomposition. The potential form minimized on that space is $m_{\kappa,\acc,0}$ by Theorem~\ref{thm:accessibility}. The index equality of Theorem~\ref{thm:separable} yields the stated count.
\end{proof}

For $j\ne0$, the bare form is
\begin{equation}\label{eq:sector-form}
 m_{\kappa,j}[a]=\int_0^R\left(
 |a'+ra/\delta^2|^2+\nu_j^2\frac{\delta^2}{r^2}|a|^2
                         -K_\eta|a|^2\right)r\dd r.
\end{equation}
The magnetic energy of high harmonics controls their current term. Put $K_+=\norm{(K_\eta)_+}_\infty$ and $R^{-2}=0$ when $R=\infty$.

\begin{corollary}[High harmonics and the full axes criterion]\label{cor:axes}
The entire direct sum of harmonics with $j\ne0$ and $\nu_j^2(1+\kappa^2/R^2)\ge K_+$ has imaginary spectrum. In particular if $(2\pi/L)^2(1+\kappa^2/R^2)\ge K_+$, the full helical spectrum lies on the two axes, and every positive real eigenvalue belongs to harmonic zero. For \eqref{eq:Gaussian-profile} under \eqref{eq:Gaussian-neutrality}, an upper bound for $K_+$ is $4d_0/(3\sqrt3\,e^{3/2})$.
\end{corollary}
\begin{proof}
Formula \eqref{eq:sector-form} gives a nonnegative bare form on the whole indicated sum. Applying the energy estimate of Theorem~\ref{thm:general-spectral} directly on this reducing subspace gives one polynomial group bound for the sum and hence a resolvent throughout $\Re\lambda\ne0$. Theorem~\ref{thm:zero} handles the complementary zero harmonic when all nonzero ones satisfy the inequality. In \eqref{eq:Gaussian-current}, the factors $(h-1)h^{-3/2}$ and $(t-1)_+e^{-t/2}$ are bounded by $2/(3\sqrt3)$ and $2e^{-3/2}$, respectively.
\end{proof}

For example, take $\kappa=m_\pm=\beta_\pm=\gamma_\pm=1$, $q_\pm=\pm1$, $c_\pm=(2\pi)^{-3/2}$, and $L=2\pi$. Then $d_0=2$ gives $K_+\le 8/(3\sqrt3\,e^{3/2})<1$. For every central value $\alpha_0>0$, only the zero harmonic of this whole-cylinder family can grow, and its growing count is exactly $N_0$.

For the whole-cylinder Gaussian branch one can also localize the exact Hamiltonian index formula. Integration by parts in \eqref{eq:sector-form} gives radial derivative energy and the potential
$(r^2-2\kappa^2)/\delta^4-K_\eta$, in addition to the nonnegative angular term. This potential is positive outside a sufficiently large disk by Proposition~\ref{prop:Gaussian-global}. Neumann bracketing of this integrated form leaves a compact-resolvent interior problem. Consequently the total bare and accessible indices are finite. The normalized accessible energy nevertheless has zero in its essential spectrum. To see this, take $a_s=s^{-1}f(r/s)$ with $f$ smooth, supported in $(1,2)$ and normalized in $L^2(r\dd r)$. Then $a_s\rightharpoonup0$, $\norm{C_e^*C_ea_s}=O(s^{-2})$, and $\norm{(\Gamma_\eta a_s,0,-\tau a_s)}\to0$. In the completed square for the accessible energy choose the particle and perpendicular field coordinates equal to $-P_{\cZ_\kappa}(\Gamma_\eta a_s,0,-\tau a_s)$, the parallel electric field zero and the potential coordinate $(I+C_e^*C_e)^{1/2}a_s$. These states have norm tending to one. The energy operator applied to the potential coordinate gives $(I+C_e^*C_e)^{-1/2}(C_e^*C_e-K_\eta)a_s\to0$. The remaining coordinates tend strongly to zero. The potential coordinates are weakly null since $\|(\Lambda-I)a_s\|\le\|C_e^*C_ea_s\|$. Compression to $\cX$ can only decrease the norm of the image, so these states give the required singular Weyl sequence.

On each fixed $j\ne0$, \eqref{eq:sector-form} gives $C_{e,j}^*C_{e,j}\ge\nu_j^2I$. The normalized free energy is coercive. The current term is compact relative to the magnetic form by local Rellich compactness and the decay of $K_\eta$. These properties pass to the accessible compression, so the exact index theorem \cite[Theorems~2.2--2.3]{LinZengIndex2022} applies to each real Fourier pair $\{j,-j\}$. Denote its accessible index by $N_{\{j,-j\}}$ and its imaginary and zero indices as in \eqref{eq:Krein-def}. Combining the pairwise index formula with Theorem~\ref{thm:zero} and Corollary~\ref{cor:axes} yields
\begin{equation}\label{eq:global-count}
 k_r+2k_c=N_0+
 \sum_{\substack{j\ge1\\ \nu_j^2<K_+}}
 \big(N_{\{j,-j\}}-2k_{i,j}^{\le0}-k_{0,j}^{\le0}\big).
\end{equation}
Only finitely many harmonics require the additional imaginary- and zero-eigenspace information. At finite radius the graph embedding into $L^2$ is compact, so bounded $K_\eta$ gives the energy gap required by \eqref{eq:HK} directly on the entire cell.

\subsection{Two Gaussian screw-pinches with imaginary spectrum}
Take $\kappa=1$, $L=\zeta=2\pi$ (one axial period per turn, hence $\nu_j=j$), $m_\pm=\beta_\pm=\gamma_\pm=1$, $q_\pm=\pm1$ and $c_\pm=5(2\pi)^{-3/2}$. Then $a_*=b_*=1$ and $d_0=10$. Two choices of the central potential illustrate why the particle invariants matter.

\begin{figure}[tb]
\centering
\includegraphics[width=.82\linewidth]{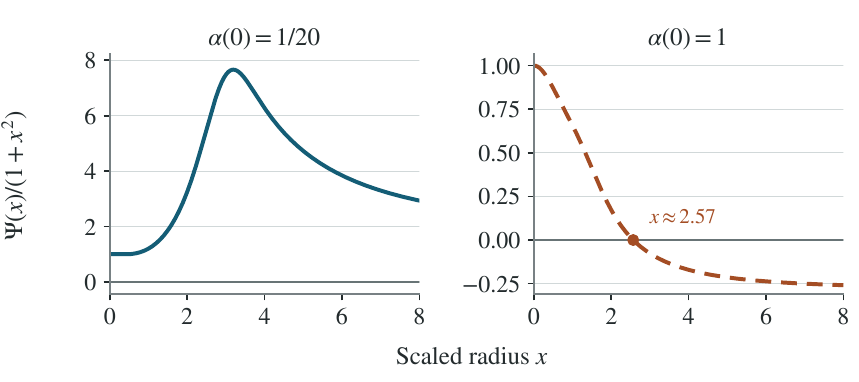}
\caption{Regular radial Jacobi fields for Theorem~\ref{thm:screw-examples}. Each panel shows $\Psi/(1+x^2)$, where $\Psi=\partial_{\alpha_0}Y$. This rescaling preserves the zeros of the regular solution $\Psi/\sqrt{1+x^2}$. The two panels have different vertical scales. For $\alpha_0=1$ the plotted zero lies at $x\simeq2.57$. Appendix~\ref{app:screw-certificate} proves the zero counts on the whole half-line.}\label{fig:screw-jacobi}
\end{figure}

\begin{theorem}[Spectral stability for two central potentials]\label{thm:screw-examples}
For both $\alpha(0)=1/20$ and $\alpha(0)=1$, the helical linearized Vlasov--Maxwell operator has spectrum contained in $i\R$ on the full potential graph space and on its Gauss/flux-constrained and accessible subspaces. For the first equilibrium the bare potential form is nonnegative. For the second it has one negative direction, while the accessible form is nonnegative.
\end{theorem}
\begin{proof}
Write $x=r$, $H=2+x^2$, $Y=\sqrt{1+x^2}\,\alpha$ and $p=x/(1+x^2)$. Let $K_*(x)=K_\eta(x)$ from \eqref{eq:Gaussian-current} for the stated parameters. After integration by parts the operator for the bare form in harmonic $j$ is
\[
 \mathsf T_j=-\partial_x^2-x^{-1}\partial_x+\frac{x^2-2}{(1+x^2)^2}
              +j^2\frac{1+x^2}{x^2}-K_*(x).
\]
At $j=0$ the regular zero-energy solution is $(1+x^2)^{-1/2}\partial_{\alpha_0}Y$. Figure~\ref{fig:screw-jacobi} illustrates the different zero counts. Indeed differentiation of \eqref{eq:Gaussian-Volterra} gives $(p(\partial_{\alpha_0}Y)')'=-pK_*\partial_{\alpha_0}Y$. Its zeros count the negative eigenvalues by Sturm oscillation on finite disks and exhaustion. Appendix~\ref{app:screw-certificate} verifies that the regular solution stays positive for the first equilibrium and crosses zero once for the second, and that a positive regular zero-energy solution exists for $|j|=1$ in both examples. The ground-state identity gives nonnegativity for $|j|=1$, and the angular term in \eqref{eq:sector-form} increases with $|j|$.

For $\alpha(0)=1$, use one conserved particle number in \eqref{eq:helical-macc}. Since $w_s=F_s$, the vector $((F_s)_s,0,0)$ belongs to $\cZ_{\kappa,0}$ and has squared norm $\mathsf M=\sum_s\iint F_s\dd v\,x\dd x$ per unit orbit measure. Its pairing with $(\Gamma_\eta a,0,-\tau a)$ is $\ip{j_\eta^0}{a}$. Hence
\[
 m_{\kappa,\acc,0}[a]\ \ge\ m_{\kappa,0}[a]
                      +\mathsf M^{-1}|\ip{j_\eta^0}{a}|^2 .
\]
The forced radial equation and its cutoff argument in Appendix~\ref{app:screw-certificate} give $m_{\kappa,0}[a]\ge Q_*^{-1}|\ip{j_\eta^0}{a}|^2$, with $3<\mathsf M<4$ and $-9<Q_*<-7$. The coefficient $\mathsf M^{-1}+Q_*^{-1}$ is positive. The accessible energy is now nonnegative on every harmonic, and thus on their orthogonal sum. Theorem~\ref{thm:helical-application} gives the full spectral conclusion.
\end{proof}

The two equilibria have the same thermal parameters, pitch and axial period. For the second, the current response gives the bare potential form one negative direction. Conservation of particle number supplies a positive term large enough to make the accessible form nonnegative. This proves the spectral conclusion without computing an eigenvalue of the phase-space generator.

\subsection*{Acknowledgments} We thank Daniel W. Crews for drawing our attention to the work of Cordonnier and collaborators and to the Liouville equation in the rigidly rotating $\theta$-pinch.

\appendix
\section{Domains, constraints, and the helical index}\label{app:framework}

\paragraph{The generator.}
In the normalized energy space \eqref{eq:normalized-space}, the skew
operator and its domain are
\begin{equation}\label{eq:Jclosed}
 \cJ z=\begin{pmatrix}
 P_+[-D(g+h)]+\Gamma_\VV E^\VV\\
 -\Gamma_\HH^*h+C_mb\\
 -\Lambda E^\VV\\
 P_-[-D(g+h)]+\Gamma_\HH E^\HH\\
 -\Gamma_\VV^*g+\Lambda\widehat a+\tau b\\
 -C_m^*E^\HH-\tau E^\VV
 \end{pmatrix},
\end{equation}
with domain
\begin{equation}\label{eq:Jdomain}
 g+h\in D(D),\quad E^\HH\in D(C_m^*),\quad b\in D(C_m),
 \quad \widehat a,E^\VV\in D(\Lambda).
\end{equation}
For two generators omit the last slot and the $C_m,\tau$ entries. Then $E^\HH$ has no derivative domain condition.

\begin{proof}[Proof of Theorem~\ref{thm:closed}]
The kinetic block is $-D$ in parity coordinates, with $g+h\in D(D)$.
The flow assumption makes it skew-adjoint. On their product domains the
field blocks $\left(\begin{smallmatrix}0&C_m\\-C_m^*&0\end{smallmatrix}\right)$
and $\left(\begin{smallmatrix}0&-\Lambda\\\Lambda&0\end{smallmatrix}\right)$
are skew-adjoint because $C_m$ is closed and $\Lambda$ is self-adjoint.
The remaining entries are bounded and skew-adjoint. Since
$(L_a-I)\widehat a=-\Lambda^{-1}(I+K)\Lambda^{-1}\widehat a$ lies in
$D(\Lambda)$, we have $D(\cJ\cL)=D(\cJ)$ and
\begin{equation}\label{eq:bounded-perturbation}
 \cA=\cJ+F,\qquad Fz=(0,0,0;0,-(I+K)\Lambda^{-1}\widehat a,0).
\end{equation}
The bounded perturbation theorem gives the group. For domain vectors
$\Lambda L_a\widehat a=Ma$, which recovers the physical field equation.
Pairing with $\cL z$ proves conservation on the domain and hence by
density on $\cE$. The bound \eqref{eq:Sigma-bound} gives a bounded
inverse change of density with $a$ retained.
\end{proof}

\paragraph{Gauss's law and fixed flux.}
Write $Gz=\diver E^\HH-Q(g+\Gamma_\VV a)\in\cT^*$ for the Gauss defect.
For $\chi\in\cT$ the domain identities are
\begin{equation}\label{eq:gradient-identities}
 DQ^*\chi=\Gamma_\HH\nabla\chi,\qquad C_m^*\nabla\chi=0,
 \qquad (Q^*\chi)_s=q_sw_s\chi.
\end{equation}
The second identity applies when the magnetic variable is present. The representing vector of $\ip{Gz}{\chi}$ is
\[
 \gamma_\chi=(-Q^*\chi,-\nabla\chi,-\Lambda^{-1}\Gamma_\VV^*Q^*\chi;0,0,0),
\]
and \eqref{eq:Jclosed} gives $\cJ\gamma_\chi=0$. Similarly, each $\zeta\in\Ker C_m$ yields a kernel vector
\[
 \varphi_\zeta=(0,0,-\Lambda^{-1}\tau\zeta;0,0,\zeta)
\]
representing $\ip{b-\tau a}{\zeta}$. Let $N$ be the closed span of these vectors. The physical homogeneous constraint space is
\begin{equation}\label{eq:Y}
 \cY=N^\perp=\{Gz=0,\ b-\tau a\in\cl{\Ran C_m^*}\}.
\end{equation}

For smooth scalar tests the transport identity follows by differentiating $q_sw_s\chi$ along the equilibrium flow. The value of $\chi$ is continuous along a reflected trajectory and $w_s$ is invariant across the reflection. Approximation in $H^1$, boundedness of $Q^*$ and $\Gamma_\HH$, and closedness of $D$ give the identity for every $\chi\in\cT$. For smooth tests with zero wall trace the tangential trace of $\nabla\chi$ vanishes and its curl is zero. Passing to the $C_m^*$ graph limit gives the second identity on the conducting domain. For two generators the $b$ constraint drops out.
More generally let $N$ be any closed subspace of $\Ker\cJ$, and let
$P$ project onto $N^\perp$. Then $P$ preserves $D(\cJ)$ and
$\cJ P=P\cJ=\cJ$. The restriction to $N^\perp$ is the compression
$\cA_P=\cJ_P(P\cL|_{N^\perp})$. Since $\cA=\cJ+F$ with $F$ bounded, the
full operator has the triangular form
\begin{equation}\label{eq:triangular}
 \cA=\begin{pmatrix}\cA_P&F_N\\0&0\end{pmatrix}
 \quad\hbox{on }P\cE\oplus N,
 \qquad F_N\hbox{ bounded}.
\end{equation}

For $\lambda\ne0$ solving the second row first gives the equivalence of
resolvents, and every generalized eigenvector has zero $N$ component.
Taking first the closed span of the Gauss and flux vectors and then $N=\Ker\cJ$ proves
\eqref{eq:nonzero-spectrum}. Exact flux representability is also preserved.
Indeed integration of the group equation gives
\begin{equation}\label{eq:exact-reconstruction}
 \int_0^tE^\HH(s)\dd s\in D(C_m^*),\qquad
 b(t)-\tau a(t)=b(0)-\tau a(0)-C_m^*\int_0^tE^\HH(s)\dd s.
\end{equation}
If $b(0)-\tau a(0)=C_m^*A^\HH(0)$, take
$A^\HH(t)=A^\HH(0)-\int_0^tE^\HH(s)\dd s$ to reconstruct an admissible potential in temporal gauge. At $\lambda\ne0$, every eigenvector satisfies
$b-\tau a=-\lambda^{-1}C_m^*E^\HH$, and induction gives the same exact-range property for generalized eigenvectors.

\paragraph{Nonzero physical normal modes.}
Consider a mode $e^{\lambda t}(f,E,B)$ of the linearized equations
\eqref{eq:linear-V}, \eqref{eq:Maxwell}--\eqref{eq:Gauss} in a retained symmetry class, with $\lambda\ne0$, $f\in\cK$, and $E,B$ in the corresponding electric and magnetic curl domains. On a conducting cylinder we take zero tangential electric trace, the maximal magnetic curl domain, and specular reflection. The particle equation is understood weakly along the reflected characteristics. On the whole cylinder we use the corresponding boundaryless domains. We assume \eqref{eq:Sigma-bound}.

The temporal-gauge mode potential is $A=-E/\lambda$. Its symmetry component $a=-E^\VV/\lambda$ lies in $L^2$, and Faraday's law gives its curl and electric graph domain. With one generator write $b=B^\VV$. The two Maxwell equations give
\[
 \begin{aligned}
 C_ea&=B^\HH, & b-\tau a&=-\lambda^{-1}C_m^*E^\HH,\\
 C_e^*B^\HH&=\lambda E^\VV+j^\VV-\tau b, &
 C_m b&=\lambda E^\HH+j^\HH .
 \end{aligned}
\]
With two generators $B^\HH=0$ by Faraday's law, $B^\VV=C_ea$ by \eqref{eq:curl-two}, and $C_e^*B^\VV=\lambda E^\VV+j^\VV$. In each case $a\in D(C_e^*C_e)$, so $\widehat a\in D(\Lambda)$.

Put $k=f-\Sigma a\in\cK$ and $\Gamma E=\Gamma_\HH E^\HH+\Gamma_\VV E^\VV$. The Noether identity gives $(\lambda+D)k=\Gamma E$ along the equilibrium characteristics. At a specular wall the factors $\mu_{s,p_j}$ agree on the incoming and reflected velocities, since $e_s$ and $p_{sj}$ are conserved at reflection. Thus the correction has the same trace from both sides. In the characteristic weak formulation this equation gives the $\cK$-valued identity
\[
 U(-t)k-k=\int_0^t U(-s)(\Gamma E-\lambda k)\,\dd s .
\]
Strong continuity gives $k\in D(D)$ on taking the difference quotient at zero. The parity components of $k$ and the fields therefore satisfy \eqref{eq:Jdomain} and the eigenvalue equation for $\cA$. This proves the normal-mode assertion on the physical domains just specified. The full-spectrum equivalence \eqref{eq:nonzero-spectrum} compares the three spaces inside $\cE$.

\begin{proof}[Proof of Theorem~\ref{thm:general-spectral}]
On $\cX$, \eqref{eq:completed-square} and energy conservation bound $\norm{E^\VV(t)}^2$ by the nonnegative initial energy $q_0$. Thus
$\norm{a(t)}\le\norm{a(0)}+|t|\sqrt{q_0}$. The normalized state norm satisfies
\begin{equation}\label{eq:norm-identity}
 \norm{z(t)}_\cE^2=q_0+\norm{a(t)}^2+\ip{Ka(t)}{a(t)}.
\end{equation}
This proves linear group growth on $\cX$. Integration of the bounded upper-right block in \eqref{eq:triangular}, now with $N=\Ker\cJ$, gives quadratic growth on the full space. If $m\ge0$, the same estimate applies directly there. The two-sided Laplace resolvent integrals exclude both open half-planes.

At finite accessible index, let $N_*=-M_\acc 1_{(-\infty,0)}(M_\acc)$, a bounded nonnegative finite-rank operator. Replace $K$ by $K-2N_*$. The Poisson operator and accessibility are unchanged, while the effective form becomes nonnegative. The auxiliary group therefore has polynomial growth and imaginary spectrum. Its generator differs from $\cA$ by a bounded finite-rank operator. Analytic Fredholm theory on each open half-plane, using invertibility sufficiently far on the real axis, gives discreteness and finite algebraic multiplicities of all nonimaginary spectrum.

On $\cX$, differentiation of the conserved form on pairs of domain vectors gives $\cA_\cX^*\cL_\cX=-\cL_\cX\cA_\cX$ on $D(\cA_\cX)$. On a nonzero root space, $\cL_\cX$ is injective, since $\Ker\cL_\cX\subset\Ker\cA_\cX$. It maps this root space at $\lambda$ into that of $\cA_\cX^*$ at $-\lambda$. Adjoint Riesz projections therefore give a root space at $-\bar\lambda$ of at least the original dimension. Reversing the roles gives equality. Real coefficients give conjugation symmetry on the complexification of the real state space. For a finite sum $E^u$ of right-half-plane root spaces, backward decay and energy conservation show $q|_{E^u}=0$. Also $E^u\cap\Ker\cL_\cX=\{0\}$. If $x_1,\ldots,x_d$ is a complex basis, choose $y_j$ with $q(x_i,y_j)=\delta_{ij}$. The energy matrix on their span is $\left(\begin{smallmatrix}0&I\\I&H\end{smallmatrix}\right)$ and has $d$ negative directions. Hence $d\le N_\acc$. Applying this to every finite selection proves \eqref{eq:unstable-bound} and finiteness of the unstable spectrum.
\end{proof}

\paragraph{The index identity with an energy gap.}
The general \emph{exact} count needs a stronger energy assumption. Suppose $\cL_\cX$ has finite negative index, finite-dimensional kernel, and a positive spectral subspace on which it is bounded below by $\delta I$ for some $\delta>0$. If $x\in\Ker\cL_\cX$, then $x=P_\cX(I-\cL)x$. The range of $I-\cL$ lies in $D(\cJ)$ by \eqref{eq:bounded-perturbation}, and $P_\cX$ preserves that domain. Hence $\Ker\cL_\cX\subset D(\cJ_\cX)$. The Lin--Zeng index theorem \cite[Theorems 2.2--2.3]{LinZengIndex2022} then gives
\begin{equation}\label{eq:HK}
 k_r+2k_c+2k_i^{\le0}+k_0^{\le0}=N_\acc.
\end{equation}
Here $E_\lambda=\bigcup_{j\ge1}\Ker(\cA_\cX-\lambda)^j$ and
\begin{equation}\label{eq:Krein-def}
 k_i^{\le0}=\sum_{\omega>0}n^{\le0}(q|_{E_{i\omega}}),\qquad
 k_0^{\le0}=n^{\le0}(q|_{E_0/\Ker\cL_\cX}).
\end{equation}
The nonpositive index retains null directions of a restricted form. At isolated finite-multiplicity nonzero imaginary eigenvalues the root-space form is nondegenerate, so $n^{\le0}$ can be replaced there by $n^-$. At embedded imaginary eigenvalues that replacement needs proof. The zero quotient in \eqref{eq:Krein-def} is essential. An odd $N_\acc$ forces a growing real eigenvalue only after the zero contribution is accounted for.

On a bounded representative cell with the intrinsic axis domain and homogeneous conducting potential data, the graph norm controls the intrinsic $H^1$ norm. Indeed invariant fields tangent to the symmetry orbits are divergence-free and
\[
 \norm{\nabla a}^2=\norm{\curl a}^2
 =\norm{C_ea}^2+\norm{\tau a}^2
\]
in the one-generator case, with the last term absent for two generators. Rellich compactness makes $\Lambda^{-1}$ compact, so $\cL=I-\text{compact}$ and the same holds after accessible compression. This verifies the stronger energy hypothesis. On an unbounded cell finite negative index alone does not supply its positive gap.

\section{Verification of the two Gaussian screw-pinches}\label{app:screw-certificate}

We prove the scalar estimates used in Theorem~\ref{thm:screw-examples}.
In the units of that theorem put $x=r$, $z=\alpha(0)$, $d=d_0=10$,
$\delta=(1+x^2)^{1/2}$, $p=x/\delta^2$, $H=2+x^2$ and $Y=\delta\alpha$.
Write $\Psi=\partial_zY$ and $W=d\,\partial_dY$, taking the derivatives
with the other parameter fixed. Differentiating \eqref{eq:Gaussian-Volterra}
gives
\begin{align}
 F=pY',\qquad &F'=pVY, & V&=d\,\delta^2H^{-3/2}e^{-\Theta/2},
           &\Theta&=Y^2/H,\nonumber\\
 G=p\Psi',\qquad &G'=-pK\Psi, & K&=V(\Theta-1),
           &&\label{eq:cert-odes}\\
 J=pW',\qquad &J'=F'-pKW. &&&\nonumber
\end{align}
Here $(Y,\Psi,W)(0)=(z,1,0)$ and the three fluxes vanish at zero.
The regular Volterra equations give differentiability in $z,d$ near the
axis, and ordinary smooth dependence continues these derivatives to
every compact interval of positive radii.

\paragraph{The zero counts.}
Write $m_j$ for the bare form in harmonic $j$ and $\mathsf T_j$ for
its operator in $L^2(x\,dx)$ with the regular axis condition. Its
differential expression is
\[
 \mathsf T_j=-\partial_x^2-x^{-1}\partial_x+
       (x^2-2)/\delta^4+j^2\delta^2/x^2-K.
\]
The regular solution of $\mathsf T_0a_0=0$ with $a_0(0)=1$ is
$a_0=\Psi/\delta$. On a disk with Dirichlet condition at the outer wall,
radial Sturm oscillation counts the negative eigenvalues by the zeros
of $a_0$ strictly before that wall. Exhaustion gives the half-line form
index since the regular compactly supported functions are a form core.
For the first harmonic set $a_1=xh$. Its regular equation is
\begin{equation}\label{eq:cert-first}
 h''+3h'/x=Qh,\qquad Q=1+(x^2-2)/\delta^4-K,\qquad h(0)=1,\quad h'(0)=0.
\end{equation}
If $h$ stays positive, the ground-state identity reads
$m_1[a]=\int_0^\infty a_1^2|(a/a_1)'|^2x\,dx$ on the regular compact
core. The boundary term at the axis vanishes since $a=O(x)$ there.
Increasing $j^2$ then gives $m_j\ge0$ for every $|j|\ge1$.

\begin{table}[htbp]
\centering\small
\caption{Strict enclosures at $R=16$ for the regular solutions of
\eqref{eq:cert-odes} and \eqref{eq:cert-first}.}\label{tab:cert-screw}
\begin{tabular}{@{}cccccc@{}}
\toprule
$z$ & $Y(R)$ & $F(R)$ & $\Psi(R)$ & $G(R)$ & $e^{-R}a_1'(R)$\\
\midrule
$1/20$ & $(334,335)$ & $(2.63,2.64)$ & $(519,520)$ & $(3.40,3.41)$ & $(0.96,0.98)$\\
$1$ & $(327,328)$ & $(2.50,2.51)$ & $(-72,-71)$ & $(-0.57,-0.56)$ & $(0.10,0.11)$\\
\bottomrule
\end{tabular}
\end{table}

Table~\ref{tab:cert-screw} gives the endpoint bounds. On $[0,16]$, $h>0$
for both equilibria and $\Psi>0$ for $z=1/20$. For $z=1$, a zero of
$\Psi$ is possible only on $[2.5625,2.578125]$, where $G<0$ throughout.
Since $\Psi(0)=1$, $\Psi(16)<0$ and $\Psi'=G/p$, there is exactly one
zero and it is simple. We next give the enclosures establishing these
claims and the quadratures used below. The supplementary source is
\texttt{verification/verify\_screw.cpp}.

\paragraph{Enclosures at the axis.}
Use $s=x^2$ up to $s_0=1/4096$. With
$N=5H^{-3/2}e^{-Y^2/(2H)}Y$, the $s$-derivatives of $F,G,J$ are
$N,N_Y\Psi,N+N_YW$, respectively. Each pair $(U,P)=(Y,F),(\Psi,G),(W,J)$
satisfies $U_s=(1+s)P/(2s)$, and $4s h_{ss}+8h_s=Qh$.
Construct degree-32 polynomials for each flux and for $h$ by matching
the source coefficients through degree 31. In subscripts denoting
coefficients of $s^k$, the recurrences for $1\le k\le32$ are
\[
 P_k=\frac{[s^{k-1}]P_s}{k},\qquad U_k=\frac{P_k+P_{k-1}}{2k},
 \qquad h_k=\frac{[s^{k-1}](Qh)}{4k(k+1)}.
\]
Set $U_{33}=P_{32}/66$ so that the companion relation holds identically.
The program encloses the real coefficients of these polynomials and checks
$\sum_k |c_k|\rho^k<2$ with $\rho=1/4$ for each $U$ and for $h$.
On the complex disk $|s|\le\rho$, this gives $|H|\ge7/4$,
$|1+s|\ge3/4$ and $|\exp(-Y^2/(2H))|<4$ on the polynomial ranges.
In particular $|N|<18$, $|N_Y|<30$, $|K|<72$ and $|Qh|<160$.
All four sources have modulus below 512. The polynomial differential
left sides have degree at most 31 and match the corresponding Taylor
coefficients. By Cauchy's estimate the absolute residual on $[0,s_0]$
is therefore at most
\[
 D=512\,(s_0/\rho)^{32}/(1-s_0/\rho).
\]

On the real box $|Y|,|\Psi|,|W|,|h|\le2$, $0\le s\le s_0$, we have
$0\le\Theta\le2$, $|N_Y|\le2$, $|N_{YY}|\le6$, $|K|\le4$,
$|K_Y|\le12$ and $|Q|\le7$. Thus $L=512$ bounds the Lipschitz constant
of the four sources in the maximum norm. Integrating a uniformly bounded
flux source to its companion costs at most $b=(s_0+s_0^2/2)/2$.
For the $h$ equation, $(s^2h_s)'=sQh/4$ gives the corresponding bound
$s_0/8$. With $bL<1$ the uniform errors for $Y,\Psi,W,h$ are at most
$E=bD/(1-bL)$. At $s_0$ the flux errors are bounded by
$s_0(LE+D)$ and the error in the $x$-derivative $h'(\sqrt{s_0})$ by
$\sqrt{s_0}(LE+D)/4$. The program checks that all three bounds are
below $2^{-289}$ and that the enlarged polynomial ranges stay strictly
inside the real box, closing the estimate by continuation. It also
checks positivity of $\Psi,h$ on the entire axis interval.

\paragraph{Enclosures from $1/64$ to $16$.}
There are eight unknowns,
$(Y,F,\Psi,G,W,J,h,h')$. Their right side $b_1$ is given by
\eqref{eq:cert-odes} and \eqref{eq:cert-first}. On an integer grid of
spacing $1/1024$ we take dyadic steps of length
$\ell\le\min(1/64,x/16)$, reaching $16$ in 1,072 steps.
At a left-endpoint box $u_0$, the program constructs a box $T$ satisfying
\[
 u_0+[0,\ell]b_1([x,x+\ell],T)\subset\operatorname{int}T .
\]
Local existence and continuation keep the solution in $T$ over the
step. Substitution of Taylor series in the differential equation gives
the coefficients through degree 19 at the left endpoint. The same
recursion over $[x,x+\ell]\times T$ bounds the twentieth derivative
divided by $20!$ at every point of the step. Taylor's formula then gives
an enclosure at the right endpoint and one on the whole step. The latter
proves the sign tests and encloses each quadrature integrand, which is
multiplied by the exact step length. The initial pieces of the
quadratures use the polynomial ranges in $s$.

The implementation uses binary interval endpoints with at least 64
significand bits. For each elementary arithmetic operation it enlarges
the rounded extrema by one representable number outward. Square roots
are bracketed by testing outward-rounded squares. For an exponential
the argument is reduced by powers of two to $|t|\le1/16$. The degree-18
Taylor polynomial there has absolute remainder at most
$2(1/16)^{19}/19!<2^{-128}$, and interval squaring recovers the original
exponential. Build instructions, floating-point requirements and the
output of a successful run accompany the source in the supplement.

\paragraph{Continuation to infinity.}
Since $F=pY'>0$ is increasing, the bounds $Y(16)>327$ and $F(16)>2.50$
give $Y(x)\ge1.1x^2$ for $x\ge16$, by the initial inequality and
$Y'(x)\ge2.50x$. Thus $\Theta\ge x^2$. Monotonicity of
$t e^{-t/2}$ and $\sqrt t\,e^{-t/2}$ for $t\ge16^2$ gives
\[
 0\le K\le10x e^{-x^2/2},\qquad 0\le pK,\ F',\ xn\le10 e^{-x^2/2},
 \qquad n=10H^{-1/2}e^{-\Theta/2}.
\]
With $R=16$, integration by parts yields
\[
 \int_R^\infty 10(1+x^2)e^{-x^2/2}\,dx
     \le 10 e^{-R^2/2}(R+2/R)=:B<2^{-167}.
\]
The endpoint checks include $|G|,|J|<10$ and $|\Psi|,|W|<10^4$.
As long as $|G|,|J|\le100$, integration of $p^{-1}=x+1/x$ gives
$|\Psi|,|W|\le1000x^2$. The total changes of each of $G,J$ on this
interval are at most $1001B<2^{-150}$ by \eqref{eq:cert-odes}.
This closes the bootstrap and keeps $G$ strictly of the same sign as
$\Psi(R)$, so no further zero occurs. The zeroth-order coefficient of
$\mathsf T_1$ is positive for $x\ge R$ since $K<1$. While $a_1>0$
we therefore have $(xa_1')'>0$. As $a_1(R),a_1'(R)>0$, the solution
stays positive to infinity. We have proved $n^-(m_0)=0$ for $z=1/20$,
$n^-(m_0)=1$ for $z=1$, and $m_j\ge0$ in both cases for $j\ne0$.

\paragraph{Removal of the negative radial direction for $z=1$.}
Here $-\partial_e F_s=F_s$ for both equilibrium densities.
The particle-number vector with kinetic components $(F_+,F_-)$ and
zero field components has squared norm $M=\int_0^\infty nx\,dx$.
Its components are constant along the equilibrium characteristics and
even under reversal of the two velocities perpendicular to the helix.
Their perpendicular current vanishes by this parity. The vector therefore
belongs to the joint kernel in harmonic zero. Its pairing with the shift produced by $a$
is $-\langle f,a\rangle$, where $f=VY/\delta$ and the pairing uses
$x\,dx$. Projection onto this vector in \eqref{eq:helical-macc}
therefore gives
\begin{equation}\label{eq:cert-number}
 m_{\mathrm{acc},0}[a]\ge m_0[a]+|\langle f,a\rangle|^2/M .
\end{equation}

Let $c=J(\infty)/G(\infty)$, which is well defined by the estimates
above, and put $u=(-W+c\Psi)/\delta$. Equations~\eqref{eq:cert-odes}
give $\mathsf T_0u=f$. The flux $p(\delta u)'=-J+cG$ tends to zero
with Gaussian decay up to a polynomial factor. Hence $\delta u$ has
a finite limit and $u=O(x^{-1})$. For a smooth cutoff $\chi_T$ equal
to one on $[0,T]$ and zero beyond $2T$, with $|\chi_T'|\lesssim T^{-1}$,
set $u_T=\chi_Tu$. Integration by parts at the regular axis gives
\[
 m_0[u_T]=\int_0^\infty \chi_T^2 u f\,x\,dx+
                 \int_0^\infty |\chi_T'|^2 u^2 x\,dx .
\]
The last term is $O(T^{-2})$, and $m_0(a,u_T)=\langle a,f\rangle$ for
all sufficiently large $T$ for every compact-core test $a$. In particular,
\[
 \lim_{T\to\infty}m_0[u_T]
   = Q_*=-\int_0^\infty F'W\,dx+c\int_0^\infty F'\Psi\,dx .
\]
The three quadratures for $M,Q_*$ have tails bounded in absolute value
by $1001B$. Adding $[-2^{-150},2^{-150}]$ separately to each finite
quadrature and to $J(R),G(R)$ gives the strict half-line bounds
\begin{align*}
 3.6402&<M<3.8736, & -1.4224&<c<-1.4223,\\
 3.1344&<\int_0^\infty F'W\,dx<3.7335, &
 2.8635&<\int_0^\infty F'\Psi\,dx<3.2492,\\
 -8.3550&<Q_*<-7.2075 .
\end{align*}
Applying Lemma~\ref{lem:rankone} to $u_T$ yields
$m_0[a]\ge |\langle f,a\rangle|^2/Q_*$. Since $Q_*<-M$,
\eqref{eq:cert-number} proves $m_{\mathrm{acc},0}\ge0$.

\bibliographystyle{amsalpha}
\bibliography{pinch_stability_strauss}

@article{Bennett1934,
  author = {Bennett, Willard H.},
  title = {Magnetically Self-Focussing Streams},
  journal = {Phys. Rev.},
  volume = {45}, number = {12}, pages = {890--897}, year = {1934},
  doi = {10.1103/PhysRev.45.890}
}

@misc{HolmerZhang2026,
  author = {Holmer, Justin and Zhang, Katherine Zhiyuan},
  title = {Equilibria for the {Vlasov--Maxwell} system related to plasma confinement},
  year = {2026}, eprint = {2608.06780}, archivePrefix = {arXiv},
  primaryClass = {math.AP}, url = {https://arxiv.org/abs/2608.06780},
  howpublished = {Preprint, \href{https://arxiv.org/abs/2608.06780}{arXiv:2608.06780}}
}

@article{Lin2022,
  author = {Lin, Zhiwu},
  title = {Linear instability of {Vlasov--Maxwell} systems revisited---a {Hamiltonian} approach},
  journal = {Kinet. Relat. Models},
  volume = {15}, number = {4}, pages = {663--679}, year = {2022},
  doi = {10.3934/krm.2021048}, eprint = {2107.12441}, archivePrefix = {arXiv}
}

@article{LinStrauss2007,
  author = {Lin, Zhiwu and Strauss, Walter A.},
  title = {Linear stability and instability of relativistic {Vlasov--Maxwell} systems},
  journal = {Comm. Pure Appl. Math.},
  volume = {60}, number = {5}, pages = {724--787}, year = {2007},
  doi = {10.1002/cpa.20158}
}

@article{LinStrauss2008,
  author = {Lin, Zhiwu and Strauss, Walter A.},
  title = {A sharp stability criterion for the {Vlasov--Maxwell} system},
  journal = {Invent. Math.},
  volume = {173}, number = {3}, pages = {497--546}, year = {2008},
  doi = {10.1007/s00222-008-0122-1}, eprint = {physics/0702023}, archivePrefix = {arXiv}
}

@book{LinZengIndex2022,
  author = {Lin, Zhiwu and Zeng, Chongchun},
  title = {Instability, index theorem, and exponential trichotomy for linear {Hamiltonian PDEs}},
  series = {Memoirs of the American Mathematical Society},
  volume = {275}, number = {1347}, publisher = {American Mathematical Society},
  address = {Providence, RI}, year = {2022},
  doi = {10.1090/memo/1347}, eprint = {1703.04016}, archivePrefix = {arXiv}
}

@article{LinZengSeparable2022,
  author = {Lin, Zhiwu and Zeng, Chongchun},
  title = {Separable {Hamiltonian PDEs} and turning point principle for stability of gaseous stars},
  journal = {Comm. Pure Appl. Math.},
  volume = {75}, number = {11}, pages = {2511--2572}, year = {2022},
  doi = {10.1002/cpa.22027}, eprint = {2005.00973}, archivePrefix = {arXiv}
}

@article{MarsdenWeinstein1982,
  author = {Marsden, Jerrold E. and Weinstein, Alan},
  title = {The {Hamiltonian} structure of the {Maxwell--Vlasov} equations},
  journal = {Phys. D},
  volume = {4}, number = {3}, pages = {394--406}, year = {1982},
  doi = {10.1016/0167-2789(82)90043-4}
}

@article{Morrison1980,
  author = {Morrison, Philip J.},
  title = {The {Maxwell--Vlasov} equations as a continuous {Hamiltonian} system},
  journal = {Phys. Lett. A},
  volume = {80}, number = {5--6}, pages = {383--386}, year = {1980},
  doi = {10.1016/0375-9601(80)90776-8}
}

@article{Morrison1981,
  author = {Weinstein, Alan and Morrison, Philip J.},
  title = {Comments on: The {Maxwell--Vlasov} equations as a continuous {Hamiltonian} system},
  journal = {Phys. Lett. A},
  volume = {86}, number = {4}, pages = {235--236}, year = {1981},
  doi = {10.1016/0375-9601(81)90496-5}
}

@incollection{GlasseyStrauss1984,
  author = {Glassey, Robert T. and Strauss, Walter A.},
  title = {Remarks on collisionless plasmas},
  booktitle = {Fluids and Plasmas: Geometry and Dynamics (Boulder, Colo., 1983)},
  series = {Contemp. Math.}, volume = {28}, publisher = {Amer. Math. Soc.},
  address = {Providence, RI}, year = {1984}, pages = {269--279},
  doi = {10.1090/conm/028/751989}
}

@article{LinStraussNonlinear2007,
  author = {Lin, Zhiwu and Strauss, Walter A.},
  title = {Nonlinear stability and instability of relativistic {Vlasov--Maxwell} systems},
  journal = {Comm. Pure Appl. Math.}, volume = {60}, number = {6},
  year = {2007}, pages = {789--837}, doi = {10.1002/cpa.20161}
}

@article{NguyenStrauss2014,
  author = {Nguyen, Toan T. and Strauss, Walter A.},
  title = {Linear stability analysis of a hot plasma in a solid torus},
  journal = {Arch. Ration. Mech. Anal.}, volume = {211}, number = {2},
  pages = {619--672}, year = {2014}, doi = {10.1007/s00205-013-0680-2}
}

@article{Antonov1961,
  author = {Antonov, V. A.},
  title = {Remarks on the problem of stability in stellar dynamics},
  journal = {Soviet Astron.}, volume = {4}, pages = {859--867}, year = {1961}
}

@article{BatchelorDavidson1975,
  author = {Batchelor, D. B. and Davidson, R. C.},
  title = {Kinetic description of linear theta-pinch equilibria},
  journal = {J. Plasma Phys.}, volume = {14}, number = {1}, pages = {77--92},
  year = {1975}, doi = {10.1017/S002237780002554X}
}

@article{SymonSeylerLewis1982,
  author = {Symon, K. R. and Seyler, C. E. and Lewis, H. R.},
  title = {Stability of {Vlasov} equilibria. {Part} 1. {General} theory},
  journal = {J. Plasma Phys.}, volume = {27}, number = {1}, pages = {13--24},
  year = {1982}, doi = {10.1017/S0022377800026349}
}

@article{LewisSeyler1982,
  author = {Lewis, H. R. and Seyler, C. E.},
  title = {Stability of {Vlasov} equilibria. {Part} 2. {One} non-ignorable co-ordinate},
  journal = {J. Plasma Phys.}, volume = {27}, number = {1}, pages = {25--35},
  year = {1982}, doi = {10.1017/S0022377800026350}
}

@article{Guo1997,
  author = {Guo, Yan},
  title = {Stable magnetic equilibria in collisionless plasmas},
  journal = {Comm. Pure Appl. Math.}, volume = {50}, number = {9},
  pages = {891--933}, year = {1997}
}

@article{Guo1999,
  author = {Guo, Yan},
  title = {Stable magnetic equilibria in a symmetric collisionless plasma},
  journal = {Comm. Math. Phys.}, volume = {200}, pages = {211--247}, year = {1999}
}

@article{MorrisonPfirsch1989,
  author = {Morrison, Philip J. and Pfirsch, Dieter},
  title = {Free-energy expressions for {Vlasov} equilibria},
  journal = {Phys. Rev. A}, volume = {40}, pages = {3898--3910},
  year = {1989}, doi = {10.1103/PhysRevA.40.3898}
}

@article{CorreaPfirsch1997,
  author = {Correa-Restrepo, Dar{\'i}o and Pfirsch, Dieter},
  title = {Negative-energy perturbations in general axisymmetric and helical {Maxwell--Vlasov} equilibria},
  journal = {Phys. Rev. E}, volume = {55}, pages = {7449--7456},
  year = {1997}, doi = {10.1103/PhysRevE.55.7449}
}

@article{Cordonnier2022,
  author = {Cordonnier, Aur{\'e}lien and Leoncini, Xavier and Dif-Pradalier, Guilhem and Garbet, Xavier},
  title = {Full self-consistent {Vlasov--Maxwell} solution},
  journal = {Phys. Rev. E}, volume = {106}, pages = {064209}, year = {2022},
  doi = {10.1103/PhysRevE.106.064209}
}

@article{Cordonnier2026,
  author = {Cordonnier, Aur{\'e}lien and Lebouazda, Yohann and Leoncini, Xavier and Dif-Pradalier, Guilhem},
  title = {Magnetically confined charged particles: {From} steep density profiles to the breaking of the adiabatic invariant},
  journal = {EPL}, volume = {155}, pages = {14002}, year = {2026},
  doi = {10.1209/0295-5075/ae6890}
}

@misc{Crews2026,
  author = {Crews, Daniel W.},
  title = {Statistical structure of canonical kinetic equilibrium with sheared flow},
  howpublished = {Preprint, \href{https://arxiv.org/abs/2607.25158}{arXiv:2607.25158}},
  year = {2026}, eprint = {2607.25158}, archivePrefix = {arXiv}
}

@article{CrewsMeierShumlak2025,
  author = {Crews, Daniel W. and Meier, Eric T. and Shumlak, Uri},
  title = {Effects of transitional orbit magnetization on transport and current in {Z} pinches},
  journal = {Phys. Plasmas}, volume = {32}, pages = {102510}, year = {2025},
  doi = {10.1063/5.0286743}
}

@article{Marcinkiewicz1939,
  author = {Marcinkiewicz, J{\'o}zef},
  title = {Sur une propri{\'e}t{\'e} de la loi de {Gau{\ss}}},
  journal = {Math. Z.}, volume = {44}, pages = {612--618}, year = {1939},
  doi = {10.1007/BF01210677}
}

@article{SchindlerPfirschWobig1973,
  author = {Schindler, K. and Pfirsch, D. and Wobig, H.},
  title = {Stability of two-dimensional collision-free plasmas},
  journal = {Plasma Phys.}, volume = {15}, number = {12},
  pages = {1165--1184}, year = {1973}, doi = {10.1088/0032-1028/15/12/001}
}

@article{Agren1986,
  author = {{\AA}gren, O.},
  title = {Variational principle and comparison theorem for the kinetic stability of two-dimensional {Vlasov} equilibria},
  journal = {Plasma Phys. Control. Fusion}, volume = {28}, number = {8},
  pages = {1055--1065}, year = {1986}
}

@article{Sharma1983,
  author = {Sharma, A. S.},
  title = {{Vlasov} stability of {Bennett} equilibrium},
  journal = {Nucl. Fusion}, volume = {23}, number = {11},
  pages = {1493--1502}, year = {1983}, doi = {10.1088/0029-5515/23/11/007}
}

@article{GuoStrauss2000,
  author = {Guo, Yan and Strauss, Walter A.},
  title = {Magnetically created instability in a collisionless plasma},
  journal = {J. Math. Pures Appl. (9)}, volume = {79}, number = {10},
  pages = {975--1009}, year = {2000}, doi = {10.1016/S0021-7824(00)01186-7}
}

@article{NguyenStrauss2013,
  author = {Nguyen, Toan T. and Strauss, Walter A.},
  title = {Stability analysis of collisionless plasmas with specularly reflecting boundary},
  journal = {SIAM J. Math. Anal.}, volume = {45}, number = {2},
  pages = {777--808}, year = {2013}, doi = {10.1137/110859695}
}

@article{Zhang2019,
  author = {Zhang, Katherine Zhiyuan},
  title = {Linear stability analysis of the relativistic {Vlasov--Maxwell} system in an axisymmetric domain},
  journal = {SIAM J. Math. Anal.}, volume = {51}, number = {6},
  pages = {4683--4723}, year = {2019}, doi = {10.1137/18M1206825}
}

@article{BenArtziHolding2017,
  author = {Ben-Artzi, Jonathan and Holding, Thomas},
  title = {Instabilities of the relativistic {Vlasov--Maxwell} system on unbounded domains},
  journal = {SIAM J. Math. Anal.}, volume = {49}, number = {5},
  pages = {4024--4063}, year = {2017}, doi = {10.1137/15M1025396}
}

@inproceedings{DickmanMorseNielson1968,
  author = {Dickman, D. O. and Morse, R. L. and Nielson, C. W.},
  title = {Numerical simulation of axisymmetric, collisionless, finite $\beta$ plasma},
  booktitle = {Proceedings of the {APS} Topical Conference on Numerical Simulation of Plasma},
  publisher = {Los Alamos Scientific Laboratory}, address = {Los Alamos, NM},
  year = {1968}, note = {Report LA-3990, contribution C2},
  url = {https://digital.library.unt.edu/ark:/67531/metadc1055154/}
}

@article{FreidbergMorse1969,
  author = {Freidberg, J. P. and Morse, R. L.},
  title = {Collisionless tearing mode instabilities in a high-$\beta$ theta pinch},
  journal = {Phys. Fluids}, volume = {12}, number = {4},
  pages = {887--894}, year = {1969}, doi = {10.1063/1.1692572}
}

@techreport{Morse1969,
  author = {Morse, Richard L.},
  title = {Equilibria of Collisionless Plasma. Part {II}},
  institution = {Los Alamos Scientific Laboratory},
  address = {Los Alamos, NM}, number = {LA-3844-MS, Part II},
  year = {1969}, url = {https://www.osti.gov/servlets/purl/4766501}
}

@article{Virtamo1980,
  author = {Virtamo, J. T.},
  title = {Locally {Maxwellian} {Vlasov} equilibria},
  journal = {Phys. Scr.}, volume = {22}, number = {2},
  pages = {134--140}, year = {1980}, doi = {10.1088/0031-8949/22/2/008}
}

@article{Meierovich1986,
  author = {Meierovich, B. E.},
  title = {Stability of a {Bennett} pinch},
  journal = {Sov. Phys. JETP}, volume = {63}, number = {5},
  pages = {966--971}, year = {1986},
  url = {https://jetp.ras.ru/cgi-bin/dn/e_063_05_0966.pdf}
}

@article{MorrisonPfirsch1990,
  author = {Morrison, Philip J. and Pfirsch, Dieter},
  title = {The free energy of {Maxwell--Vlasov} equilibria},
  journal = {Phys. Fluids B}, volume = {2}, number = {6},
  pages = {1105--1113}, year = {1990}, doi = {10.1063/1.859246}
}

@article{KiesslingLebowitz1994,
  author = {Kiessling, Michael K.-H. and Lebowitz, Joel L.},
  title = {Dissipative stationary plasmas: {Kinetic} modeling, {Bennett}'s pinch and generalizations},
  journal = {Phys. Plasmas}, volume = {1}, number = {6},
  pages = {1841--1849}, year = {1994}, doi = {10.1063/1.870639}
}

@article{ChannonCoppins2001,
  author = {Channon, S. W. and Coppins, M.},
  title = {{Z}-pinch {Vlasov-fluid} equilibria including sheared-axial flow},
  journal = {J. Plasma Phys.}, volume = {66}, number = {5},
  pages = {337--347}, year = {2001}, doi = {10.1017/S0022377801001507}
}

@article{KomarovFadeev1962,
  author = {Komarov, N. N. and Fadeev, V. M.},
  title = {Plasma in a self-consistent magnetic field},
  journal = {Sov. Phys. JETP}, volume = {14}, number = {2},
  pages = {378--381}, year = {1962},
  note = {Russian original, Zh. Eksp. Teor. Fiz. 41 (1961), 528--533},
  url = {https://jetp.ras.ru/cgi-bin/dn/e_014_02_0378.pdf}
}

@article{Pfirsch1962,
  author = {Pfirsch, D.},
  title = {{Mikroinstabilit{\"a}ten vom Spiegeltyp in inhomogenen Plasmen}},
  journal = {Z. Naturforsch. A}, volume = {17}, number = {10},
  pages = {861--870}, year = {1962}, doi = {10.1515/zna-1962-1004}
}

@inproceedings{Schindler1966,
  author = {Schindler, K.},
  title = {A variational principle for one-dimensional plasmas},
  booktitle = {Proceedings of the Seventh International Conference on Phenomena in Ionized Gases},
  publisher = {Gradevinska Knjiga}, address = {Belgrade}, volume = {II},
  pages = {736--739}, year = {1966},
  url = {https://juser.fz-juelich.de/record/830307}
}

@article{BellCamposPintoEtAl2026,
  author = {Bell, Dominik and Campos Pinto, Martin and Possanner, Stefan and Sonnendr{\"u}cker, Eric},
  title = {The linearized {Vlasov--Maxwell} system as a {Hamiltonian} system},
  journal = {J. Comput. Phys.}, volume = {555}, pages = {114765},
  year = {2026}, doi = {10.1016/j.jcp.2026.114765}
}

@article{MondtDavidson1980,
  author = {Mondt, J. P. and Davidson, R. C.},
  title = {Hybrid-kinetic stability properties of high-beta plasmas},
  journal = {J. Plasma Phys.}, volume = {23}, number = {1},
  pages = {157--201}, year = {1980}, doi = {10.1017/S0022377800022224}
}

@article{HolmMarsdenRatiuWeinstein1985,
  author = {Holm, Darryl D. and Marsden, Jerrold E. and Ratiu, Tudor and Weinstein, Alan},
  title = {Nonlinear stability of fluid and plasma equilibria},
  journal = {Phys. Rep.}, volume = {123}, number = {1--2},
  pages = {1--116}, year = {1985}, doi = {10.1016/0370-1573(85)90028-6}
}

@article{KruskalOberman1958,
  author = {Kruskal, M. D. and Oberman, C. R.},
  title = {On the stability of plasma in static equilibrium},
  journal = {Phys. Fluids}, volume = {1}, number = {4},
  pages = {275--280}, year = {1958}, doi = {10.1063/1.1705885}
}
\end{document}